\documentclass[10pt,a4paper]{article}
\usepackage{amsfonts}
\usepackage{amssymb}
\usepackage{mathrsfs}
\usepackage{amsthm}
\usepackage{setspace}
\usepackage{color}

\usepackage{fancyhdr}
\newtheorem{theo}{Theorem}[section]
\newtheorem{prop}[theo]{Proposition}
\newtheorem{cor}[theo]{Corollary}
\newtheorem{lemm}[theo]{Lemma}

\theoremstyle{definition}
\newtheorem*{rem}{Remark}
\newtheorem{defi}[theo]{Definition}

\newenvironment{lproof}{\emph{Proof of Lemma.}}{ \qed \par}

\newcommand{\be}{\begin{eqnarray*}}
\newcommand{\ee}{\end{eqnarray*}}
\newcommand{\beqa}{\begin{eqnarray}}
\newcommand{\eeqa}{\end{eqnarray}}
\newcommand{\ba}{\begin{array}}
\newcommand{\ea}{\end{array}}
\newcommand{\onab}{\overrightarrow{\nabla}}
\newcommand{\mc}{\mathcal}
\newcommand{\mf}{\mathfrak}
\newcommand{\ora}{\overrightarrow}
\newcommand{\rP}{\mathsf{P}}

\newcommand{\mbb}{\mathbb}
\newcommand{\wt}{\widetilde}
\newcommand{\wh}{\widehat}

\newcommand{\iso}{\textrm{iso}}

\newcommand\CR{\begin{picture}(46,12)%
\put(5,3){\line(1,0){16}}%
\put(25,3){\line(1,0){16}}%
\put(3,3){\makebox(0,0){$\times$}}%
\put(23,3){\makebox(0,0){$\circ$}}%
\put(43,3){\makebox(0,0){$\times$}}%
\end{picture}}

\newcommand\Bpalgebra{\begin{picture}(76,12)\put(5,3){\line(1,0){16}}%
\put(25,3){\line(1,0){6}}%
\put(51,3){\line(-1,0){6}}\put(55,4){\line(1,0){17}}%
\put(54,2){\line(1,0){18}}\put(39,3){\makebox(0,0){\dots}}%
\put(63,2){\makebox(0,0){$>$}}%
\put(3,3){\makebox(0,0){$\circ$}}%
\put(23,3){\makebox(0,0){$\circ$}}%
\put(53,3){\makebox(0,0){$\circ$}}%
\put(73,3){\makebox(0,0){$\times$}}%
\end{picture}}

\newcommand\Dpalgebra{\begin{picture}(76,12)\put(5,3){\line(1,0){16}}%
\put(25,3){\line(1,0){6}}%
\put(51,3){\line(-1,0){6}}%
\put(55,4){\line(2,1){17}}%
\put(54,2){\line(2,-1){18}}%
\put(39,3){\makebox(0,0){\dots}}%
\put(3,3){\makebox(0,0){$\circ$}}%
\put(23,3){\makebox(0,0){$\circ$}}%
\put(53,3){\makebox(0,0){$\circ$}}%
\put(74,12){\makebox(0,0){$\circ$}}%
\put(73,-7){\makebox(0,0){$\times$}}%
\end{picture}}

\newcommand\Dpalgebrajaw{\begin{picture}(76,12)\put(5,3){\line(1,0){16}}%
\put(25,3){\line(1,0){6}}%
\put(51,3){\line(-1,0){6}}%
\put(55,4){\line(2,1){17}}%
\put(54,2){\line(2,-1){18}}%
\put(39,3){\makebox(0,0){\dots}}%
\put(3,3){\makebox(0,0){$\circ$}}%
\put(23,3){\makebox(0,0){$\circ$}}%
\put(53,3){\makebox(0,0){$\circ$}}%
\put(74,12){\makebox(0,0){$\times$}}%
\put(73,-7){\makebox(0,0){$\times$}}%
\end{picture}}

\newcommand\Balgebra{\begin{picture}(76,12)\put(5,3){\line(1,0){16}}%
\put(25,3){\line(1,0){6}}%
\put(51,3){\line(-1,0){6}}\put(55,4){\line(1,0){17}}%
\put(54,2){\line(1,0){18}}\put(39,3){\makebox(0,0){\dots}}%
\put(63,2){\makebox(0,0){$>$}}%
\put(3,3){\makebox(0,0){$\circ$}}\put(23,3){\makebox(0,0){$\circ$}}%
\put(53,3){\makebox(0,0){$\circ$}}\put(73,3){\makebox(0,0){$\circ$}}%
\end{picture}}

\newcommand\Dalgebra{\begin{picture}(76,12)\put(5,3){\line(1,0){16}}%
\put(25,3){\line(1,0){6}}%
\put(51,3){\line(-1,0){6}}%
\put(55,4){\line(2,1){17}}%
\put(54,2){\line(2,-1){18}}%
\put(39,3){\makebox(0,0){\dots}}%
\put(3,3){\makebox(0,0){$\circ$}}\put(23,3){\makebox(0,0){$\circ$}}%
\put(53,3){\makebox(0,0){$\circ$}}%
\put(74,12){\makebox(0,0){$\circ$}}%
\put(74,-7){\makebox(0,0){$\circ$}}%
\end{picture}}

\newcommand\trione{\begin{picture}(49,12)\put(24,3){\line(-1,0){18}}%
\put(28,4){\line(2,1){17}}%
\put(27,2){\line(2,-1){18}}%
\put(4,3){\makebox(0,0){$\circ$}}%
\put(26,3){\makebox(0,0){$\circ$}}%
\put(47,12){\makebox(0,0){$\circ$}}%
\put(46,-7){\makebox(0,0){$\times$}}%
\end{picture}}

\newcommand\tritwo{\begin{picture}(49,12)\put(24,3){\line(-1,0){18}}%
\put(28,4){\line(2,1){17}}%
\put(27,2){\line(2,-1){18}}%
\put(5,2){\makebox(0,0){$\times$}}%
\put(26,3){\makebox(0,0){$\circ$}}%
\put(47,12){\makebox(0,0){$\circ$}}%
\put(47,-7){\makebox(0,0){$\circ$}}%
\end{picture}}

\begin{document}

\title{Free $n$-distributions: holonomy, sub-Riemannian structures, Fefferman constructions and dual distributions.}
\author{Stuart Armstrong}
\date{2007}
\maketitle

\begin{abstract}
This paper analyses the parabolic geometries generated by a free $n$-distribution in the tangent space of a manifold. It shows that certain holonomy reductions of the associated normal Tractor connections, imply preferred connections with special properties, along with Riemannian or sub-Riemannian structures on the manifold. It constructs examples of these holonomy reductions in the simplest cases. The main results, however, lie in the free $3$-distributions. In these cases, there are normal Fefferman constructions over CR and Lagrangian contact structures corresponding to holonomy reductions to $SO(4,2)$ and $SO(3,3)$, respectively. There is also a fascinating construction of a `dual' distribution when the holonomy reduces to $G_2'$.
\end{abstract}

\tableofcontents

\newpage

\section{Introduction}

On a manifold $M$, let $H \subset TM$ be a distribution of rank $n$. Then there is a well defined map $\mc{L}: H \wedge H \to TM/H$. For $X,Y$ sections of $H$, it is given by the quotiented Lie bracket $X\wedge Y \to [X,Y]/H$. Then $H$ is a \emph{free} $n$-distribution if $\mc{L}$ is an isomorphism. The moniker ``free'' comes from the fact that there are no relations between sections of $H$ that would cause $\mc{L}$ to fail injectivity.

This condition immediately implies that $TM/H$ is of rank $n(n-1)/2$, thus that $M$ is of dimension $m = n(n+1)/2$. Bryant \cite{bryskew} has studied the case of $n=3, m=3$, a free $3$-distribution in a $6$-manifold, but the general case remains little studied.

Fortunately, these structures lead themselves to be treated with the general tools of Cartan connections on parabolic geometries (\cite{TCPG} and \cite{CartEquiv}). The homogeneous model is provided by the set of maximal isotropic planes in $\mathbb{R}^{n+1,n}$. The group of transformations is $G = SO(n+1,n)$ while the stabiliser of a point is $P = GL(n) \rtimes \mathbb{R}^n \rtimes \wedge^2 \mathbb{R}^n$. Its Lie algebra is $\mf{p}$ which has nilradical $\mathbb{R}^n \rtimes \wedge^2 \mathbb{R}^n$. These are precisely the two-step free nilpotent Lie algebras, with the Lie bracket from $\mathbb{R}^n \otimes \mathbb{R}^n$ to $\wedge^2 \mathbb{R}^n$ being given by taking the wedge. The fact this nil-radical is free is a consequence of the freeness of the $n$-distribution.

We do not introduce any extra information, or make any choices by taking the Cartan connection, as the normal Cartan connection for a free $n$-distribution is determined entirely by $H$ (\cite{two}).

The most natural restrictions to put on the holonomy of a connection with structure algebra $\mf{so}(n+1,n)$ is to require that it preserves a subundle in the natural representation bundle of that algebra -- the standard Tractor bundle $\mc{T}$. This condition is analysed; it turns out it implies a class of preferred connections on $M$, which preserve certain structures on the manifold, making them an example of sub-Riemannian manifolds.

If the rank of the preserved bundle $V \subset \mc{T}$ is $n$, there is a unique preferred connection $\nabla$ defined by it, that has properties analogous to the Einstein condition in conformal and projective geometry. If $V$ is further non-degenerate, there is a well-defined metric on the manifold as a whole. In that case, it is an Einstein involution \cite{meein}.

Other issues worth looking into in any new geometries is how the structures restrict to sub-manifolds; this is analysed in the next section.

There is even a decomposition/twisted product result, similar to the Einstein product result in conformal geometry (\cite{felcon} and \cite{mecon}) which applies to certain very restrictive holonomy algebras. In this case, there are explicit constructions of manifolds with these properties, leaving hope that manifolds with the weaker properties mentioned above will also exist.

The main results of this paper will be gleaned in the $n=3, m=6$ case. The free $3$-distribution has a Fefferman construction into the conformal structure \cite{bryskew}. We will show this Fefferman construction is normal for both Tractor connections, meaning that we have many known examples of holonomy reductions \cite{mecon}.

Here, the normal Cartan connection is torsion-free, and the results of the previous section can be applied to show that holonomy reductions to $SU(2,2) \cong Spin(4,2)_0$ and $SL(4,\mbb{R}) \cong Spin(3,3)_0$ do exist, and arise from their own Fefferman constructions -- over integrable CR manifolds and integrable Lagrangian contact structures, respectively. Here, normality of the underlying Cartan connections is equivalent with normality of that generated by the free $3$-distribution.

The other interesting situation for a free $3$-distribution is that of a reduction of the holonomy of the Tractor connection $\onab$ to $G_2'$. This does not arise from any Fefferman construction, but has a fascinating geometry. On an open dense set of the manifold, there is a canonical Weyl structure $\nabla$. This determines a splitting of $T = T_{-2} \oplus H$, where $H$ is the canonical free $3$-distribution. Then $H' = T_{-2}$ is also a free $3$-distribution, and the normal Tractor connection it generates is isomorphic with $\onab$. Iterating this procedure generates $H$ again; thus $H$ and $H'$ are in some sense `dual' distributions.

\section{The geometry of free $n$-distributions}

\subsection{Homogeneous model} \label{hom:mod}

The homogeneous model for a free $n$-distribution is the space of isotropic $n$-planes inside $\mathbb{R}^{(n+1,n)}$. Let $\mf{g} = \mf{so}(n+1,n)$, and put the metric on $\mathbb{R}^{(n+1,n)}$ in the form:
\be
\left( \begin{array}{ccc} 0&0&Id_n \\ 0&1&0 \\ Id_n &0&0 \end{array} \right),
\ee
where $Id_n$ is the identity matrix on $\mathbb{R}^n$. The algebra $\mf{g}$ is then spanned by elements of the form:
\be
\left( \begin{array}{ccc} A&v&B \\ w&0&v^t \\ C&w^t&-A^t \end{array} \right),
\ee
where $A \in \mf{gl}(n)$, $B^t = - B$ and $C^t = -C$. Then the isotropic plane
\be
V = \{(a_1, a_2, \ldots, a_n, 0, 0, \ldots, 0) | a_j \in \mathbb{R}\}
\ee
is preserved by the subalgebra $\mf{p} \subset \mf{g}$ spanned by elements of the form
\be
\left( \begin{array}{ccc} A&v&B \\ 0&0&v^t \\ 0&0&-A^t \end{array} \right),
\ee
This algebra is isomorphic to $\mf{gl}(n) \oplus (\mathbb{R}^n)^* \oplus (\wedge^2 \mathbb{R}^n)^*$, with the natural algebraic structure.

The homogeneous manifold is $M = G/P$, where $G = SO(n+1,n)_0$ and $P \subset G$ is the Lie group with Lie algebra $\mf{p}$. $M$ is of dimension $(2n+1)(2n)/2 - (n^2 + n + n(n-1)/2) = n(n+1)/2$. There is a subspace $\mf{g}_{(-1)}$ of $\mf{g}$, consisting of those elements with $C = 0$. Right action on $G$ generalises this to a distribution $\widehat{H} \subset TG$. This distribution is preserved by $P$, and hence the map $G \to M$ maps $\widehat{H}$ to a distribution $H \subset TM$. This distribution is evidently of rank $n$.
 
Now consider the differential bracket of right-invariant vector fields which are sections of $\widehat{H}$; this matches up with the Lie bracket on $\mf{g}_{(-1)}$, thus $[\widehat{H}, \widehat{H}]$ spans all of $TG$. Consequently, $[H,H]$ must span all of $TM$, making $H$ free by dimensional considerations.

To have an explicit model, let $x_j$, $1 \leq j \leq n$ and $x_{kl}$, $1 \leq k < l \leq n$ be local coordinates on $M$. Then we may define vector fields $U_{kl}$ via:
\be
U_{kl} &=& \frac{\partial}{\partial x_{kl}}.
\ee
Extend this definition by requiring $U_{lk} = - U_{kl},$ $U_{kk} = 0$. Then we complete the frame with the vector fields
\be
X_{j} &=& \frac{\partial}{\partial x_{j}} + \frac{1}{2} \sum_{p=1}^{n} x_{p} U_{pj}.
\ee
These vector fields follow the commutator relations:
\be
[X_{p}, X_{j} ] &=& U_{pj}, \\[0.cm]
[U_{k_1 l_1}, U_{k_2 l_2}] &=& 0, \\[0.cm]
[U_{kl}, X_j] &=& 0 .
\ee
In this model, the distribution $H$ is simply the span of $X_j$.

\subsection{Cartan connection}

Given a semi-simple Lie algebra $\mathfrak{g}$ with Killing form $(-,-)$, a subalgebra $\mf{p} \subset \mf{g}$ is said to be \emph{parabolic} iff $\mf{p}^{\perp}$ is the \emph{nilradical} of $\mf{p}$, i.e.~its maximal nilpotent ideal. This gives (\cite{paradef}, details are also available in the author's thesis \cite{methesis}) a filtration of $\mf{g}$:
\be
\begin{array}{ccccccccccccc}
\mf{g}_{(-k)}& \supset & \mf{g}_{(k-1)} & \supset &\ldots& \supset &\mf{g}_{(0)}& \supset &\mf{g}_{(-1)}& \supset & \ldots & \supset & \mf{g}_{(k)}, \\
| | &&&&&&| | && | | && \\
\mf{g} &&&&&&\mf{p}&& \mf{p}^{\perp}&&
\end{array}
\ee
such that $\{ \mf{g}_{(j)} , \mf{g}_{(l)} \} \subset \mf{g}_{(j+l)}$. The associated graded algebra is $ {\textrm{gr}}(\mf{g}) = \bigoplus_{-k}^{k} \mf{g}_j$, where $\mf{g}_j = \mf{g}_{(j)}/\mf{g}_{(j+1)}$. By results from \cite{paradef}, ${\textrm{gr}}(\mf{g})$ is isomorphic to $\mf{g}$. Furthermore, there is a unique element $\epsilon_0$ in $\mf{g}_0$ such that $\{ \epsilon_0, \xi \} = j \xi$ for all $\xi \in \mf{g}_j$. The isomorphisms ${\textrm{gr}}(\mf{g}) \cong \mf{g}$ compatible with the filtration are then given by a choice of lift $\epsilon$ of $\epsilon_0$ with respect to the exact sequence
\beqa \label{p:sequence}
0\to \mf{p}^{\perp} \to \mf{p} \to \mf{g}_0 \to 0.
\eeqa
This means that the gradings of $\mf{g}$ compatible with the filtration form an affine space modelled on the nilradical $\mf{p}^{\perp}$. Define $P$ and $G_0$ as the subgroups of $G$ that preserve the filtration and the grading, respectively. It is easy to see that their Lie algebras are $\mf{p}$ and $\mf{g}_0$, and that the inclusion $G_0 \subset G$ is non-canonical.

\begin{defi}[Cartan connection]
A Cartan connection on a manifold $M$ for the parabolic subalgebra $\mf{p} \subset \mf{g}$ is given by a principal $P$-bundle
\be
\mc{P} \to M,
\ee
and a one form $\omega \in \Omega^1(\mc{P}, \mf{g})$ with values in the Lie algebra $\mf{g}$ such that:
\begin{enumerate}
\item $\omega$ is equivariant under the $P$-action ($P$ acting by $Ad$ on $\mathfrak{g}$),
\item $\omega(\sigma_A) = A$, where $\sigma_A$ is the fundamental vector field of $A \in \mathfrak{p}$,
\item $\omega_u: T \mathcal{P}_u \to \mathfrak{g}$ is a linear isomorphism for all $u \in \mathcal{P}$.
\end{enumerate}
\end{defi}

The inclusion $P \subset G$ generates a bundle inclusion $\mc{P} \subset \mc{G}$ and $\omega$ is then the pull back of a unique $G$-equivariant connection form on $\mc{G}$ which we will also designate by $\omega$. Since $\omega$ takes values in $\mf{g}$, it generates a standard connection on any vector bundle associated to $\mc{G}$. This connection is called the \emph{Tractor connection} and will be designated by $\onab$.

In the case we are looking at, $\mf{g} = \mf{so}(n+1,n)$ and $\mf{p} = \mf{gl}(n) \oplus (\mathbb{R}^n)^* \oplus (\wedge^2 \mathbb{R}^n)^*$. $\mf{g}_0$ is simply $\mf{gl}(n)$, and the nilradical of $\mf{p}$ is
\be
\mf{p}^{\perp} = (\mathbb{R}^n)^* \oplus (\wedge^2 \mathbb{R}^n)^*.
\ee
And all the two-s1tep free-nilpotent algebras are precisely of this form. In terms of the notation for parabolic subalgebras introduced \cite{CartEquiv}, this is:
\Bpalgebra.

Define the Lie algebra bundle
\be
\mc{A} = \mc{P} \times_P \mf{g}.
\ee
This has a natural filtration
\be
\mc{A} =\mc{A}_{(-2)} \supset \mc{A}_{(-1)} \supset \mc{A}_{(0)} \supset \mc{A}_{(1)} \supset \mc{A}_{(2)}.
\ee
Paper \cite{TCPG} demonstrates that the tangent space $T$ of $M$ is equal to the quotient bundle $\mc{A} / \mc{A}_{(0)}$. The killing form gives an isomorphism
\be
\mc{A}_{(-1)} \cong (\mc{A} / \mc{A}_{(0)})^* = T^*.
\ee
Hence there is a well defined inclusion $T^* \subset \mc{A}$, and a well defined projection $\mc{A} \to T$.

\begin{defi}[Weyl structure]
A \emph{Weyl structure} on $(M, \mc{P}, \omega)$ is a $P$-equivarient function $\eta: \mc{P} \to \mathfrak{p}$ that is always a lift of the grading element $\epsilon_0$, as in equation (\ref{p:sequence}).
\end{defi}
A Weyl structure gives a splitting of $\mf{g}$, and consequently allows a decomposition of both the Lie algebra bundle and the Cartan connection:
\be
\mc{A} &=& \mc{A}_{-2} \oplus \mc{A}_{-1} \oplus \mc{A}_{0} \oplus \mc{A}_{1} \oplus \mc{A}_{2} \\
\omega &=& \omega_{-2} + \omega_{-1} + \omega_0 + \omega_1 + \omega_2.
\ee
Dividing out by the action of $\mf{p}^{\perp}$ gives a quotient map $\mf{p} \to \mf{g}_0$ and hence a bundle map $\mc{P} \to \mc{G}_0$. There is a unique $G_0$-equivarient one-form on $\mc{G}_0$ of which $\omega_0$ is the pull-back; we will call it $\omega_0$ as well. This is a principal connection on $\mf{G}_0$; since $G_0$ acts on $\mf{g}_{-2} + \mf{g}_{-1} \cong \mf{g}/\mf{p}$, then $T = \mc{A} / \mc{A}_{(0)}$ is an associated bundle to $\mc{G}_0$. This implies that $\omega_0$ generates an affine connection $\nabla$ on the tangent bundle.

These $\nabla$'s are called preferred connections; they are in one-to-one correspondence with Weyl structures and hence to compatible splittings of $\mc{A}$. They consequently form an affine space modelled on $\mc{A}_{(-1)} = T^*$; the relation between two preferred connections $\nabla$ and $\widehat{\nabla}$ is given explicitly by the one-form $\Upsilon$ on $M$ such that
\be
\nabla_X Y = \widehat{\nabla}_X Y + \{\{X, \Upsilon \}, Y\}_{-},
\ee
with $\{,\}$ the Lie bracket on $\mc{A}$. The splitting of $\mc{A}$ gives further splittings $T = T_{-2} \oplus T_{-1}$ and $T^*_{1} \oplus T^*_{2}$. The bundles $T_{-1}$ and $T^*_{2} = (T_{-1})^{\perp}$ are defined independently of the Weyl structure, since they are preserved by the action of $T^* \subset \mc{A}$.

Given a preferred $\nabla$, the Tractor connection on $\mc{A}$ and any associated bundles is given by:
\be
\onab_X v = \nabla_X v + X \cdot v + \rP(X) \cdot v,
\ee
where $\rP$ is the rho-tensor, a section of $T^* \otimes T^*$, and $\cdot$ is the action of $T$ and $T^*$ on $\mc{A}$ given by the Lie bracket.

\begin{defi}[Curvature]
The curvature of the Cartan connection is defined to be the two-form $\kappa = d \omega + \{ \omega , \omega \} \in \Omega^2(\mc{P}, \mf{g})$. It is easy to see that $\kappa$ vanishes on vertical vectors, and is $P$-equivariant; consequently dividing out by the action of $P$, $\kappa$ may seen as an element of $\Omega^2(M, \mathcal{A})$; in this setting, it is the curvature of the Tractor connection $\onab$. Finally, the inclusion $T^* \subset \mc{A}$ implies that $\kappa$ is equivalent to a $P$-equivariant function from $\mc{P}$ to $\wedge^2 \mf{p}^{\perp} \otimes \mf{g}$. We shall use the designation $\kappa$ interchangeably for these three equivalent definitions, though it is generally the third one we shall be using.
\end{defi}

Given a grading on $\mf{g}$, there is a decomposition of any tensor product $\otimes^c \mf{g} = \sum \mf{g}_{j_1, j_2, \ldots j_c} $ where $\mf{g}_{j_1, j_2, \ldots j_c} = \mf{g}_{j_1} \otimes \mf{g}_{j_2} \otimes \ldots \otimes \mf{g}_{j_c}$.

The homogeneity of $\mf{g}_{j_1, j_2, \ldots j_c}$ is defined to be the sum $j_1 + j_2 + \ldots + j_c$. Any element of $\eta$ of $\otimes^c \mf{g}$ can be decomposed into homogeneous elements $\eta_{j_1, j_2, \ldots, j_c}$. The \emph{minimal homogeneity} of $\eta$ is defined to be the lowest homogeneity among the non-zero $\eta_{j_1, j_2, \ldots, j_c}$.

Homogeneity is not preserved by the action of $P$; however since $\mf{p}$ consists of elements of homogeneity zero and above, the minimal homogeneity of any element is preserved by the action of $P$. Since $\kappa$ is a map to $\wedge^2 \mf{p}^{\perp} \otimes \mf{g}$, the following definition makes sense:
\begin{defi}[Regularity]
A Cartan connection is regular iff the minimal homogeneity of $\kappa$ is greater than zero.
\end{defi}

There are well defined Lie algebra differentials $\partial: \wedge^c \mathfrak{p}^{\perp} \otimes \mf{g} \to \wedge^{c+1} \mathfrak{p}^{\perp} \otimes \mf{g}$ and codifferentials $\partial^*: \wedge^c \mathfrak{p}^{\perp} \otimes \mf{g} \to \wedge^{c-1} \mathfrak{p}^{\perp} \otimes \mf{g}$. In terms of decomposable elements, the co-differential is given by
\be
\partial^* (u_1 \wedge \ldots \wedge u_c) \otimes v &=& \sum_{j \neq k} u_1 \wedge \ldots \wedge \{u_j, u_k \} \wedge \ldots \wedge u_c \otimes v + \\
&& \sum_{j} u_1 \wedge \ldots  \hat{u_j} \ldots \wedge u_c \otimes \{u_j, v\}.
\ee
\begin{defi}[Normality] \label{normal:def}
A Cartan connection is normal iff $\partial^* \kappa = 0$.
\end{defi}
If we let $(Z_l)$ be a frame for $T$ and $(Z^l)$ a dual frame for $T^*$, this condition is given, in terms of $\kappa$ an element of $\Omega^2(M, \mc{A})$, as
\be
(\partial^* \kappa)(X) = \sum_{l} \{ Z^l, \kappa_{(Z_l, X)} \} - \frac{1}{2} \kappa_{(\{Z^l,X \}_{-}, Z_l)},
\ee
for all $X \in \Gamma(T)$.

Paper \cite{CartEquiv} demonstrates that a regular, normal Cartan connection for these Lie groups is determined entirely by the distribution $T_{-1}$. We shall call this distribution $H$, as before. Since $\omega$ is regular, the algebraic bracket $\{ ,\}$ matches up with graded version of the Lie bracket $[,]$ of vector fields; in other words, if $X$ and $Y$ are sections of $H$,
\be
[X,Y] - \{X,Y\} \in \Gamma(H),
\ee
implying that $H$ is maximally non-integrable. Indeed, any maximally non-integrable $H$ of correct dimension and co-dimension determines a unique normal Cartan connection as above.

\begin{defi}
We shall call $(M,H)$ a manifold with a free $n$-distribution. It is always of dimension $m = n(n+1)/2$.
\end{defi}

\subsection{Harmonic curvature} \label{har:cur}

A Cartan connection is said to be \emph{torsion-free} if the curvature $\kappa$ seen as a function $\mc{P} \to \wedge^2 \mf{p}^{\perp} \otimes \mf{g}$ is actually a function $\mc{P} \to \wedge^2 \mf{p}^{\perp} \otimes \mf{p}$. If $\kappa$ is instead seen as a section of $\wedge^2 T^* \otimes \mc{A}$, torsion-freeness implies that is a section of $\wedge^2 T^* \otimes \mc{A}_{(1)}$.

Since $\partial^* \kappa = 0$, $\kappa$ must map into Ker$(\partial^*)$. This has a projection onto the homology component $H_2(\mf{p}^{\perp}, \mf{g}) =$ Im$(\partial^*)/$Ker$(\partial^*)$. The composition of $\kappa$ with this projection is $\kappa_H$, the \emph{harmonic curvature}. The Bianci identity for normal Cartan connections imply that it is a complete obstruction to integrability (\cite{CartEquiv}). Indeed, paper \cite{BGG} demonstrates that $\kappa_H$ determines $\kappa$ entirely.

Now, Kostant's solution of the Bott-Borel-Weil theorem \cite{Kostant} allows us to algorithmically calculate $H_2(\mf{p}^{\perp}, \mf{g})$. The $n = 1$ case is trivial, the $n=2$ and $n=3$ cases have $H_2(\mf{p}^{\perp}, \mf{g})$ contained inside
\be
(\mf{g}_{1} \wedge \mf{g}_{2}) \otimes \mf{g}_{0}.
\ee
The harmonic curvature must lie inside this component, which is of homogeneity three. Both the Bianci identity for normal Cartan connections \cite{CartEquiv}, and the construction of the full curvature from the harmonic curvature \cite{BGG} imply that the other component of the curvature must have higher homogeneity. Since the torsion components have maximal homogeneity three (for $(\wedge^2 \mf{g}_2) \otimes \mf{g}_{-1}$), these two geometries are torsion-free.

For $n \geq 4$, the harmonic curvature is contained inside
\be
(\mf{g}_{1} \wedge \mf{g}_{2}) \otimes \mf{g}_{-2},
\ee
and thus these geometries are never torsion-free (unless they are flat). This will be most evident when we look at preserved substructures (see Section \ref{max:pref}); while the torsion-free geometries do not require extra conditions for these substructures to correspond to structures on immersed submanifolds, we will require extra conditions on the curvature $\kappa$ in the general case for this to be true.

\subsection{The Tractor bundle} \label{trac:bun}
The standard Tractor bundle $\mc{T}$ is bundle on which we will be doing most of our calculations. It is defined to be the bundle associated to $\mathbb{P}$ under the standard representation of $\mf{so}(n+1,n)$:
\be
\mc{P} \times_P \mathbb{R}^{(n+1,n)}.
\ee

This makes $\mc{T}$ into a rank $2n+1$ bundle. A choice of preferred connection $\nabla$ reduces the structure group of $\mc{T}$ to $\mf{gl}(n)$. Under this reduction,
\be
\mc{T} = H \oplus \mathbb{R} \oplus H^*.
\ee
Changing the the choice of $\nabla$ by a one-form $\Upsilon$ changes this splitting as:
\beqa \label{splitting:change}
\left( \begin{array}{c} v \\ \tau \\ X \end{array} \right) \to \left( \begin{array}{c} v + \tau \Upsilon_1 - \{\Upsilon_2, X \} - \frac{1}{2} (\Upsilon_1(X)) \Upsilon_1 \\ \tau -\Upsilon_1(X) \\ X \end{array} \right).
\eeqa
This demonstrates that the inclusions $H^* \subset \mathbb{R} \oplus H^* \subset \mc{T}$ are well defined, as are the projections $\mc{T} \to H \oplus \mathbb{R} \to H$. The projection that we shall be using the most often is $\pi^2: \mc{T} \to H$.

The metric $h$ on $\mc{T}$, of signature $(n+1, 1)$ is given in this splitting by:
\be
h \ \big( \left( \begin{array}{c} v \\ \tau \\ X \end{array} \right) , \left( \begin{array}{c} w \\ \nu \\ Y \end{array} \right) \big) = \frac{1}{2} (w(Y) + v(X) + \tau \nu),
\ee
while the Tractor derivative in the direction of any $Z \in \Gamma(T)$ is $\onab_Z = Z + \nabla_Z + \rP(Z)$, or, more explicitly,
\be \label{tractor:con}
\onab_Z \left( \begin{array}{c} v \\ \tau \\ X \end{array} \right) = \left( \begin{array}{c} \nabla_Z v + \tau \rP(Z)_1 - \{\rP(Z)_2 , X\} \\ \nabla_Z \tau - v(Z_{-1}) - \rP(Z)_1(X) \\ \nabla_Z X +  \tau Z_{-1} + \{ Z_{-2}, v \} \end{array} \right).
\ee
\begin{rem}
Note that the $Z$ term in $Z + \nabla_Z + \rP(Z)$ implies there cannot exist a section of $H^*$ or $\mathbb{R} \oplus H^*$ that is mapped into the same bundle by $\onab$, at any point. The converse of that is that if $V$ is any bundle preserved by $\onab$, then its intersection with $H^*$ and $\mathbb{R} \oplus H^*$ will be minimal on an open, dense set. As we will generally be dealing with such preserved $V$'s and as our results are local, we will generally assume that this intersection is minimal, by restricting implicitly to the open dense subset where it is true.
\end{rem}

\section{Preserved subundles} \label{pers:struc}
We shall be looking at the various implications of preserved subundles of the Tractor bundle $\mc{T}$. We shall always be assuming that $\onab$ is normal, unless explicitly stated otherwise, but most of these results do not need the normality condition.

Recall that for any bundle $A \subset B$, the bundle $A^{\perp} \subset B^*$ is defined to be the maximal bundle such that the contraction $A^{\perp} \llcorner A$ is always zero. For any metric $g$ on $A$, we define $\iso^g(A)$ to be $A \cap g(A^{\perp})$ -- the maximal isotropic subspace of $A$.

Now, let $V$ be a \emph{generic} subundle of $\mc{T}$ of rank $\leq n$. By generic we mean that the projection $\pi^2: \mc{T} \to H$ maps $V$ injectively to a subundle $A \subseteq H$, on a open dense subset of $M$. The remark at the end of Section \ref{trac:bun} implies that any bundle preserved by $\onab$ is generic. The metric $h$, restricted to $V$ and then projected to $A$, gives a (possibly very degenerate or null) metric $g$ on $A$.

A choice of preferred connection $\nabla$ gives a splitting of $\mc{T}$, and hence a map $V \to H^*$. Since $A$ and $V$ are isomorphic, this gives a map $\mu^{\nabla}: A \to H^*$. Finally, dividing out by the action of $A^{\perp}$ gives a map $g^{\nabla} : A \to A^*$.

\begin{defi}[$V$-preferred connections]
We say that the connection $\nabla$ is $V$-preferred if $g^{\nabla} = g$. We say the connection $\nabla$ is strongly $V$-preferred if it is preferred and $\mu^{\nabla}(V) \cap A^{\perp} = 0$. The difference between preferred and strongly preferred is specious if $g$ (hence $V$) is non-degenerate. For degenerate $V$, a $V$-preferred $\nabla$ has a $\mu^{\nabla}$ that maps isotropic elements of $A$ to sections of $A^{\perp}$, while a strongly $V$-preferred $\nabla$ will have a $\mu^{\nabla}$ that is zero on isotropic elements of $A$.
\end{defi}

From now on, we will drop the superscript from $\mu^{\nabla}$, referring to it as simply $\mu$ unless we want to emphasis the dependence. It is good to have a $\mu$ that is a metric on all of $T$, not just on $A$; we can extend it as follows:

If $\nabla$ is strongly $V$-preferred, we can extend $\mu$ to a map $H \to H^*$ simply by defining it to be zero on $(\mu(A))^{\perp} \subset H$. Since $\mu(A) \cap A^{\perp} = 0$, $A + (\mu(A))^{\perp} = H$ and $\mu$ is zero on $A \cap (\mu(A))^{\perp}$, making this well defined. This extended $\mu$ is still symmetric.

If $\nabla$ is just $V$-preferred, we may also extend $\mu$, but we have to make a choice of a bundle $F$ transverse to $A + (\mu(A))^{\perp} \subset H$. Then we extend $\mu$ by defining it to be zero on $(\mu(A))^{\perp}$, requiring that $\mu(F) \subset F^{\perp}$, and requiring that $\mu$ remain symmetric. The actual choice of $F$ will not be important.

Finally, notice that the metric $\mu$ extends to $T_{-2}$ by the isomorphism $T_{-2} =\{H,H\} \cong H \wedge H$. Defining $T_{-2}$ and $H$ to be perpendicular extends the definition of $\mu$ to a metric on all of $T$.

\begin{theo} \label{V:pref}
For any generic $V$ of rank $\leq n$, there exists preferred and strongly preferred connections $\nabla$ for $V$. The $V$-preferred connections form an affine space modelled on $A^{\perp} \oplus \{ H, A^{\perp} \}$; the strongly $V$-preferred form an affine space modelled on $A^{\perp} \oplus \{\left(\iso^{\mu}(A)\right)^{\perp}, A^{\perp}\}$. Furthermore, the bundle $B = \{A,A\}$, is independent of the choice of $V$-preferred connection.
 
If $\onab$ preserves $V$, these $V$-preferred connections $\nabla$ have the following properties:
\begin{enumerate}
\item $A$ is preserved by $\nabla$ along directions in $C = H \oplus B + \{\mu(A)^{\perp},\mu(A)^{\perp}\} \subset T$,
\item $\nabla \mu$ is zero on $A \otimes A$,
\item $\rP_{11} = -\mu + \eta$ with $\eta$ a section of $H^*\otimes A^{\perp}$,
\item $\rP_{22} = -\mu + \eta'$, with $\eta'$ a section of $T_2^* \otimes \{A^{\perp}, H^* \} $,
\item $\rP_{21}$ and $\rP_{12}$ are sections of $T^*_2 \otimes A^{\perp}$ and $T^*_1 \otimes \{A^{\perp}, H^* \}$, respectively.
\end{enumerate}
If $\nabla$ is actually strongly $V$-preferred, then $\eta'$ is a section of $T^*_2 \otimes \{A^{\perp}, \iso^g(A)^{\perp} \}$ and $\rP_{12}$ a section of $T^*_1 \otimes \{A^{\perp}, \iso^g(A)^{\perp} \}$.
\end{theo}

The properties of this theorem derive from the following proposition:
\begin{prop}
Let $V\subset \mc{T}$ be a generic bundle of rank $\leq n$. Then if $\nabla$ determines a splitting such that $V$ has no $\mathbb{R}$ component and $\mu$ is symmetric on $A$, then $\nabla$ is a $V$-preferred connection. If $V \cap (H, 0, A^{\perp}) = 0$, it is strongly $V$-preferred. And if $V$ is preserved by $\onab$, $\nabla$ has the properties enumerated in Theorem \ref{V:pref}. 
\end{prop}
\begin{proof}
Having no $\mathbb{R}$ component is understood to mean that given the splitting $H \oplus \mathbb{R} \oplus H^*$, projection onto the central component maps $V$ to zero. In this case, the map $\mu: A \to H^*$ is defined by mapping the $H$ component of $V$ to its $H^*$ component. If $\check{X} = (X,0, \mu(X))$ and $\check{Y} = (Y,0, \mu(Y))$ are sections of $V$, the metric $h$ on these elements is given by:
\be
g(X,Y) = h(\check{X}, \check{Y}) = \frac{1}{2} (\mu(X,Y) + \mu(Y,X)) = \mu(X,Y),
\ee
making $\nabla$ $V$-preferred. $\mu|_{A \otimes A}$ remains constant under the action of a section of $A^{\perp} \oplus \{ H, A^{\perp} \}$.

For the strongly $V$-preferred connections, $V \cap (H, 0, A^{\perp}) = 0$ iff $\mu(A) \cap A^{\perp}$ = 0. Strongly $V$-preferred connections are related by the action of $A^{\perp} \oplus \{\left(\iso^{\mu}(A)\right)^{\perp}, A^{\perp}\}$, the subspace of $A^{\perp} \oplus \{ H, A^{\perp} \}$ that preserves the above relation.

Now assume that $\onab$ preserves $V$. If $Z$ is a section of $C = H \oplus B + \{\mu(A)^{\perp},\mu(A)^{\perp}\} \subset T$, then $Z \cdot (0,0,\mu(A)) \subset (A, \mathbb{R}, 0)$. Since
\be
\onab_Z (X,0, \mu(X)) = (\nabla_Z X + \{Z_{-2}, \mu(X) \}, \ldots) \in \Gamma(V),
\ee
$\nabla_Z$ must map sections of $A$ to sections of $A$.

Now let $Z$ be any section of $T$. Pick $\check{X} = (X, 0, \mu(X))$ and $\check{Y} = (Y,0, \mu(Y))$ such that $h(\check{X}, \check{Y}) = \mu(X,Y)$ is a constant. Then since $\onab$ preserves $h$:
\be
h(\onab_Z \check{X}, \check{Y}) + h(\check{X}, \onab_Z \check{Y}) = 0
\ee
But this is:
\be
\mu(Y, \nabla_Z X + \{Z_{-2} ,\mu(X)\}) + \mu(X, \nabla_Z Y +\{ Z_{-2}, \mu(Y)\}).
\ee
The terms $\mu(Y, \{Z_{-2} ,\mu(X)\}) + \mu(X, \{Z_{-2}, \mu(Y)\})$ must vanish if $Z_{-2}$ is decomposable under the identification $T_{-2} \cong \wedge^2 H$ -- hence it must vanish for all $Z$, giving
\be
0 = \mu(Y, \nabla_Z X) + \mu(X, \nabla_Z Y),
\ee
implying that $\nabla \mu$ is zero on $A \otimes A$.

Now take the derivative of any section $\check{X}$ of $V$:
\be
\onab_Z \check{X} = (\nabla_Z X + \{ Z_{-2} \cdot \mu(X)\}, -\mu(X,Z_{-1}) - \rP_{11}(Z,X) - \rP_{21}(Z,X) , \nabla_Z \mu(X) - \{ \rP(Z)_{2}  X\}).
\ee
Since this must be a section of $V$, the central term must vanish. Consequently $\rP_{21}$ vanishes on $A$, so is a section of $T_{2}^* \otimes A^{\perp}$, and $\rP_{11} = -\mu + \eta$ with $\eta$ a section of $H^*\otimes A^{\perp}$. Also since this is a section of $V$,
\be
\mu(\nabla_Z X + \{ Z_{-2} , \mu(X)\}) = \nabla_Z \mu(X) - \{ \rP(Z)_{2} , X\}).
\ee
Then since the $\mu(\nabla_Z X) - \nabla_Z \mu(X) = (\nabla_Z \mu)(X)$ is a section of $A^{\perp}$, so must be $\mu (\{ Z_{-2} , \mu(X)\}) + \{ \rP(Z)_{2} , X\}$. This implies that $\rP_{22} + \mu$ and $\rP_{12}$ must, under the bracket action, map sections of $A$ to sections of $A^{\perp}$. This means that $\rP_{22} + \mu$ is a section of $T_2^* \otimes \{A^{\perp}, H^* \} $ and a section of $T^*_1 \otimes \{A^{\perp}, H^* \}$.

If $\nabla$ is moreover strongly $V$-preferred, then any $X$ that is a section of $\iso^g(A)$, must have $\mu(X) = 0$, and $\mu(A) \cap A^{\perp} = 0$, implying that $\{\rP(Z)_{22} + \mu(Z) ,X \} = \{ \rP(Z)_{12},  X\} = 0$. Consequently $\rP_{22} + \mu$ and $\rP_{12}$ are sections of $T^*_2 \otimes \{A^{\perp}, \iso^g(A)^{\perp} \}$ and $T^*_1 \otimes \{A^{\perp}, \iso^g(A)^{\perp} \}$, respectively.
\end{proof}
However, we have not yet shown the existence of preferred connections $\nabla$ with the required properties. This done in the next lemma:
\begin{lemm}
For any generic bundle $V$ of rank $r<n$ such that $\pi^2: V \to A \subset H$ is bijective, there locally exists a preferred connection $\nabla$ such that $V$ has no $\mathbb{R}$ component, and $\mu$ is symmetric on $A$. Moreover there exists such $\nabla$ so that $V \cap (0,\mathbb{R}, A^{\perp}) = 0$.
\end{lemm}
\begin{lproof}
Pick any frame $(v_j)$ for $V$, and any preferred connection $\nabla$. Then, assuming without loss of generality that $v_1$ has an $\mathbb{R}$ factor, we may construct another frame $(v_j')$ where $v_1'= (X_1, 1, \omega_1)$ is a scaling of $v_1$ and the other $v_j' = (X_j,0, \omega_j)$ have no $\mathbb{R}$ factors.

Then pick a section $\alpha$ of $H^*$ such that $\alpha(X_j) = \delta_{1j}$. Changing the preferred connection $\nabla$ by the action of $\alpha$ gives a splitting where $v_1'= (X_1, 1 - \alpha(X_1), \omega_1') = (X_1, 0 , \omega_1')$ and $v_j' = (X_j, -\alpha(X_j), \omega_j') = (X_j,0,\omega_j')$. So $V$ has no central component.

To make $\mu$ symmetric, pick an orthogonal frame $(w_j)$ of $V$, with $w_j = (Y_j, 0, \nu_j)$. If we can ensure that $\nu_j(Y_k) = 0$ whenever $j \neq k$, then $\mu$ must be symmetric on the frame $(Y_j)$ and hence on the whole of $A \otimes A$. We will proceed by induction; assume that $\nu_j(Y_k) = 0$ for $j \neq k$, $j \leq l$. Since the $w_j$ are orthogonal, $\nu_{l+1}(Y_{k}) = 0$ for $k < l$. We may now change the splitting by the action of the component
\be
\Upsilon = \sum_{k = l+2}^r \nu_{l+1}(Y_k) \{Y^*_{l+1}, Y^*_{k} \}
\ee
where $Y_p^*$ is any section of $H^*$ such that $Y_p^*(Y_q) = \delta_{pq}$. The action of $\Upsilon$ is trivial on all $Y_k$ with $k < j$, so this change does not change the relations $\nu_j(Y_k) = 0$ for $j \neq k$, $j \leq l$. However, in the new splitting, $\nu_{l+1}(Y_k) = 0$ for $k \neq l+1$.

Since this process works for $l = 0$ as well, by induction there exists a $\nabla$ giving a splitting with no $\mathbb{R}$ component such that $\mu$ is symmetric on $A$.

Now, assume that we have a $V$-preferred connection $\nabla$ as above, and pick a new frame $(v_j)$ such that $(v_j)_{j \leq p}$ is a frame for $\iso^h(V)$. The splitting we are looking for is one where $v_j = (X_j,0,0)$ for $j\leq p$. This result can also be proved inductively. We shall be using sections of $\{H^*, A^{\perp} \}$ to change the splitting, which will not affect the fact that $\nabla$ is $V$-preferred.

Assume that we have a splitting where $v_j = (X_j, 0,0)$ for $j \leq l$ where $0 \leq l < p$. Then $v_{l+1} = (X_{l+1}, 0, \tau)$. Since $v_{l+1}$ is orthogonal to all of $V$, $\tau$ must be a section of $A^{\perp}$. Choose a section $\xi$ of $H^*$, defined so that $\xi(X_{j}) = \delta_{jl}$. Then changing the splitting by the action of $\{ \xi, \tau\}$ gives a new splitting where $v_j = (X_j,0,0)$ for $j \leq l + 1$.

Since the construction holds for $l = 0$ as well, by induction there exists a strongly $V$-preferred $\nabla$.

\end{lproof}

\begin{theo} \label{max:pref}
If $V$ is generic of rank $n$, the conditions of Theorem \ref{V:pref} simplify considerably. In these cases, there exist a unique $V$-preferred connection $\nabla$ (which is automatically strongly $V$-preferred).

And if $\onab$ preserves $V$:
\begin{itemize}
\item $\nabla \mu = 0$,
\item $\rP_{21} = \rP_{12} = 0$,
\item $\rP_{11} = - \mu$ on $H$, and $\rP_{22} = -\mu$ on $T_{-2}$,
\item hence $\nabla \rP = 0$.
\end{itemize}
\end{theo}
If $V$ is a definite subspace of $\mc{T}$, then this is in fact an Einstein involution, so these properties are expected -- see paper \cite{meein}.

By similar arguments to above, if $V$ is preserved only in the directions \emph{along} $H$, then we have an example of sub-Riemannian Geometry, for the metric $\mu$ on $H$.

\subsection{Sub-Tractor bundles}
Here we attempt to grapple with the issue of sub-structures contained within the total structure -- specifically, of free $n$-distributions and their Cartan connections on a distribution $E \subset T$ of lower rank. Since these Cartan connections have torsion in the general case, the integrability of $E$ must be addressed separately.

We want conditions for $\onab$ to descend to a Tractor connection on a subdistribution $E \subset T$. Assume there is a bundle $A$ and a preferred connection $\nabla$ such that $\nabla$ preserves $A$ along $E = A + \{A,A\}$. We will work in the splitting given by $\nabla$, and we want $\onab$ to descend to a Tractor connection on $E$. There are two versions of this: a strong and a weak condition.

\begin{defi}[Strong condition]
$\onab$ preserves $E \oplus (A \otimes A^*) \oplus E^*\subset \mf{g}(\mc{T})$ along $E$ for some choice of $A^* \subset H^*$ (which then determines $E^*$). This implies that $\nabla$ preserves $E^*$ along $E$.
\end{defi}

\begin{defi}[Weak condition]
$\onab$ preserves $\mf{E}$ along $E$ for some choice of $A^* \subset H^*$ (which then determines $E^*$), where
\be
E \subset \mf{E} \subset E \oplus (A \otimes A^*) \oplus E^*\subset \mf{g}(\mc{T}),
\ee
for some choice of $A^* \subset H^*$. This implies that $\nabla$ preserves $E^* \cap \mf{E}$ along $E$.
\end{defi}
The strong condition describes a sub-tractor bundle $\mc{A}_E$ contained in the algebra bundle $\mc{A}$. The weak condition describes a preserved subundle of such a sub-Tractor bundle.

\begin{lemm}
For the strong and weak conditions, $E$ is integrable if and only if $\kappa(Z_1, Z_2)_j$, $j\leq1$ is a section of $E$ whenever $Z_1$ and $Z_2$ are sections of $E$.
\end{lemm}
\begin{proof}
By the properties above, $\rP(E)$ preserve $\mf{E}$. Hence if the ``torsion-terms'' of $\kappa(Z_1,Z_2)$ are sections of $E$, then the actual torsion of $\nabla$ on $E\wedge E$ is a section of $E$. Since $\nabla$ preserves $E$ along $E$, this implies that $E$ is integrable.
\end{proof}

\begin{rem}[Inclusion]
If we have the strong condition, then define $\mc{T}_E = A \oplus \mathbb{R} \oplus A^*$ and we have an inclusion
\be
\mc{T}_E \subset \mc{T}.
\ee
This bundle must now be preserved by $\onab$ along $E$. This inclusion is unchanged by changing $\nabla$ by any section $\Upsilon$ of $E^* = A^* \oplus \{A^*, A^*\}$.
\end{rem}

Our work on preserved bundles in Section \ref{pers:struc} gives us a useful class of sub-Tractor bundles:
\begin{prop}
Assume that we have a subundle $V \subset \mc{T}$ of rank $n$, preserved by $\onab$, with $V_1 \subset V$ also preserved. Then if $A = \pi^2(V_1)$ and $E = A \oplus \{A, A \}$, we have a weak sub-Tractor connection on $V_1$. If $\mu$ is non-degenerate on $V_1$, we have the strong condition, with a well defined inclusion $\mc{T}_{E} \subset \mc{T}$.
\end{prop}
\begin{proof}
Pick the unique $V$-preferred connection $\nabla$; since this is also $V_1$-preferred, it preserves $A$ along $E$. Define $A^* = \mu(A) \subset H^*$; since $\nabla \mu = 0$ this is preserved by $\nabla$ along $E$ as well. Let $Z$ be any section of $E$; consequently, $\rP(Z)$ is a section of $E_1^*$, and $\mu(E)$ is preserved along $E$. Thus starting from $E \subset \mc{A}$, repeated differentiation along $E$ generates the algebra
\be
E \oplus (A \otimes \mu(A)) \oplus \mu(E),
\ee
and $\onab$ preserves that bundle along $E$. If $\mu$ is non-degenerate, we have the strong condition, $\mu$ uniquely defined by the uniqueness of $\nabla$. Moreover the inclusion
\be
\mc{T}_E = A \oplus \mathbb{R} \oplus \mu(A) \subset H \mathbb{R} \oplus H^* = \mc{T},
\ee
is well defined.
\end{proof}

A slight weakening of the previous conditions give us:
\begin{prop}
Assume we have a preserved $V$ of rank $r \leq n$. If $\nabla_Z \mu = 0$ for $Z$ any section of $E$, we have a weak sub-Tractor bundle. If also $\mu$ is of maximal rank, then we have the strong condition.
\end{prop}
\begin{proof}
Same as previous.
\end{proof}

Finally, there is a consideration of normality. When is this sub-Tractor bundle normal, assuming $E$ is integrable? If $Z$ is a section of $E$, then:
\be
0= (\partial^* \kappa)(Z) &=& \sum_{Z_j} \{\kappa(Z_j, Z), Z^j\} - \frac{1}{2} \kappa_{(Z^j, Z)_{-}, Z_j} \\
&=& \sum_{Z'_j} \{\kappa(Z'_j, Z), (Z^j)'\} - \frac{1}{2} \kappa_{((Z^j)', Z)_{-}, Z'_j} \\
&& + \sum_{Z''_j} \{\kappa(Z''_j, Z), (Z^j)''\} - \frac{1}{2} \kappa_{((Z^j)'', Z)_{-}, Z''_j},
\ee
where $(Z'_j)$ is a frame for $E$, and $(Z''_j)$ a frame for $(E^*)^{\perp}$. The normality condition is then that $\sum_{Z''_j} \{\kappa(Z''_j, Z), (Z^j)''\} - \frac{1}{2} \kappa_{((Z^j)'', Z)_{-}, Z''_j}$ has a trivial action on $\mf{E}^{\perp}$ ($\mc{A}_E$ for the strong condition). There is no reason to suppose this is true in general, though Corollary \ref{sub:decomp} in the next section demonstrates a particular case where normality does descend to $E$.

\section{Twisted product}
In this section we present a decomposition and twisted product construction for certain types of free $n$-distributions, dependent on the holonomy and curvature of the Tractor connections. We will first need to introduce some terminology and definitions:
\begin{defi}
A connection $\onab$ \emph{fixes} a vector bundle $V$, if there is a frame $(v_j)$ of $V$ such that $\onab v_j = 0$.
\end{defi}
A co-isotropic bundle $V \subset \mc{T}$ is a bundle such that $W = h(V^{\perp})$ is contained in $V$. A minimal co-isotropic bundle is one of rank $n+1$ (which implies that $W$ of rank $n$ is a maximal isotropic bundle).

If $(M_1, H_1)$ and $(M_2, H_2)$ are two free distributions, let $N = H_1 \bigotimes H_2$ be the total space of the tensor product of $H_1$ and $H_2$ over $M_1 \times M_2$ (in other words there is a surjective map $H \to M_1 \times M_2$ and the vertical fibre of $TN$ at the point $(x,y) \in M_1 \times M_2$ is $(H_1)_x \times (H_2)_y$). The point of this construction is the following theorem:
\begin{theo} \label{comb:theo}
If $(M_1, H_1)$ and $(M_2, H_2)$ are two manifolds carrying free distributions of ranks $n_1$ and $n_2$ such that their Tractor connections $\onab^1$ and $\onab^2$ each fix minimal co-isotropic bundles $V_1$ and $V_2$, then there exists a well-defined free distribution on the total space of the tensor product $N = H_1 \bigotimes H_2$, and a well-defined normal Cartan connection $\onab^3$ for this structure -- and also fixes a minimal co-isotropic bundle $V_3$, with an isomorphism $W_3 \cong W_1 \oplus W_2$. We call this the twisted product of $(M_1, \nabla_1)$ and $(M_2, \nabla_2)$.
\end{theo}
If $\onab^1$ fixes $V_1$, choose the unique $W_1$-preferred connection $\nabla^1$ defined by the bundle $W_1 = h(V_1^{\perp})$ (see Theorem \ref{max:pref}). In this splitting, $\rP = 0$. Moreover, since $\onab^1$ fixes $V_1$, it must also fix its isotropic subspace, i.e. $W_1$. Since the projection $W_1 \to H$ is bijective, this means that $\nabla^1$ is actually a flat connection. Let $(X_j)$ be a local frame of $H$ such that $\nabla^1 X_j = 0$.

Moreover, $\onab^1$ must also fix another section of $V$. In the splitting defined by $\nabla^1$, this must be of the form $(A,1,0)$ (since it is orthogonal to all the elements $(X_j, 0,0)$). By the formula for $\onab^11$ of Equation (\ref{tractor:con}), this implies that $\nabla_{X_j}^1 A = - X_j$ and $\nabla_{U} A = 0$ for $U$ a section of $(TM_1)_{-2}$; in other words, $\nabla^1 A = Id_H$.

Now $A = x^j X_j$ for some functions $x_j$. The above formulas thus implies that $X_j \cdot x^k = - \delta^k_j$, and, furthermore, that $U \cdot x^k = 0$ for all sections $U$ of $(TM_1)_{-2}$.

Similar conclusion hold for $V_2$, generating a flat preferred connection $\nabla^2$, a frame $(Y_k)$ of $H_2$ and functions $(y^k)$ such that $Y_j \cdot y^k =\delta^k_j$.

Now, if $\mc{V}$ is the vertical subundle of $TN$, there is an exact sequence:
\be
\begin{array}{ccccccccccc} \label{exact:sequence}
0 &\to &\mc{V} &\to& TN &\to &TM_1& \times &TM_2 &\to& 0. \\
&&&&&&\cup&&\cup&& \\
&&&&&&H_1&&H_2&&
\end{array}
\ee
The preferred connections $\nabla_1$ and $\nabla_2$ give a connection on $N \to M_1 \times M_2$ seen as a vector bundle and hence a splitting of the above sequence, given by a map $\sigma: TM_1 \oplus TM_2 \to TN$. We want to adjust this splitting so that the inclusion $H_1 \oplus H_2 \to TN$ defines a maximally non-integrable subundle of $H_3$.

To do so, define sections $\widetilde{X}^j$ and $\widetilde{Y}_k$ of $TN$ by
\be
\widetilde{X}_j &=& \sigma(X_j) - \frac{1}{2} \sum_k y^k X_j \otimes Y^k \\
\widetilde{Y}_k &=& \sigma(Y_k) + \frac{1}{2} \sum_j x^j X_j \otimes Y^k.
\ee
Define $H_3$ to be the span of these elements. The projection $TN \to TM_1 \oplus TM_2$ maps $H_3$ to $H_1 \oplus H_2$. Moreover, since $\nabla_1$ and $\nabla_2$ respect the algebraic brackets on $M_1$ and $M_2$, we have the following commutator relations:
\be
[\widetilde{X}_j , \widetilde{X}_k ] &=& \sigma([X_j, X_k]) \\[0.cm]
[\widetilde{Y}_j , \widetilde{Y}_k ] &=& \sigma([Y_j, Y_k]) \\[0.cm]
[\widetilde{X}_j , \widetilde{Y}_k ] &=& X_j \otimes Y_k.
\ee
The new splitting $\sigma'$ of the sequence (\ref{exact:sequence}) is given my $\sigma'(X_j) = \widetilde{X}_j$, $\sigma'(Y_k) = \widetilde{Y}_k$ and $\sigma' = \sigma$ on $(TM_1)_{-2}$ and $(TM_2)_{-2}$.

We now define $\nabla^3$ by requiring it to be zero on all of $\widetilde{X}_j$ , $\widetilde{Y}^k$ and all of their Lie brackets. Given $\nabla^3$ and $\sigma'$, the Tractor bundle $\mc{T}_3$ is defined as $H_3 \oplus \mathbb{R} \oplus H^*_3$, and the Tractor connection as $\onab^3_Z = Z + \nabla^3_Z$. To show it is normal, we need the following proposition:
\begin{prop}
The curvature of $\onab^3$ is the direct sum of those of $\onab^1$ and $\onab^2$, and hence $\onab^3$ is normal. On a open, dense set, the local holonomy of $\onab^3$ is the sum of those of $\onab^1$ and $\onab^2$. Moreover, $\onab^3$ also fixes a maximal co-isotropic bundle.
\end{prop}
\begin{proof}
Both $\nabla^1$ and $\nabla^2$ are flat, and have vanishing $\rP$. Hence the only terms in the curvatures of $\onab^1$ and $\onab^2$ are the torsion terms of $\nabla^1$ and $\nabla^2$. We need to show that $\nabla^3$ has exactly the same torsion terms -- which, since it is flat, is equivalent with demonstrating that the Lie bracket of its preserved sections is the same.

First, it is easy to see that the splitting $\sigma$ preserves the Lie bracket on $TM_1$ and on $TM_2$. The splitting $\sigma'$ adds extra terms to $TM_1$, but all the extra terms are vertical vectors multiplied by a function whose derivative vanishes along $TM_1$ and along $(TM_2)_{-2}$. This implies that $\sigma'$ preserves the Lie bracket on $TM_1$ and does not introduce any extra torsion between $TM_1$ and $(TM_2)_{-2}$. The same argument shows that there are no extra torsion terms on $TM_2$ or between $TM_2$ and $(TM_1)_{-2}$. The Lie bracket between horizontal sections and vertical sections of $TN$ is also trivial, as it is on vertical sections of $TN$. Since we have defined the algebraic bracket in such a way that it matches the differential one on sections $\widetilde{X}_j$ and $\widetilde{Y}_k$, this demonstrates we have no extra torsion terms.

Since the co-differential $\partial^*$ is $C^{\infty}(N)$-linear, the normality of $\onab^1$ and $\onab^2$ imply the normality of $\onab^3$.

Similarly, the infinitesimal holonomy of $\onab^3$ -- the span of the image in $\mc{A}$ of the iterated derivatives of the curvature of $\onab^3$ -- is the sum of the infinitesimal holonomies of $\onab^1$ and $\onab^2$. And on an open, dense set, the infinitesimal holonomy of any connection matches up to its local holonomy.

The maximal co-isotropic bundle $V \subset \mc{T}_3$ preserved by $\onab^3$ is that spanned by $(\widetilde{X}_j,0,0)$ , $(\widetilde{Y}_k, 0,0)$ and $(\sigma'(A_1) + \sigma'(A_2), 1,0)$, where $A_j$ is the section of $H_j$ defined by $\nabla^j A_j = Id_{H_j}$.
\end{proof}

This twisted product has a corollary, a decomposition result:

\begin{cor} \label{sub:decomp}
Given any manifold $(M,H)$ with a free $n$-distribution, with normal Tractor connection $\onab$ that fixes a minimal co-isotropic bundle $V$, let $W_1$ be any preserved subspace of $W = h(V^{\perp})$, and $A_1 = \pi(W_1)$ and $E_1 = A_1 \oplus \{ A_1, A_1 \}$ in the splitting defined by the $W$-preferred connection $\nabla$. Then if $\kappa(E_1,E_1)_{-} \subset E_1$, $E_1$ is integrable. If, moreover, there is a preserved $W_2$ such that $W = W_1 \oplus W_2$, with the same properties, and $\kappa(E_j, T) = \kappa(E_j, E_j)$, for $j=1,2$ then the (weak) sub-Tractor connections on $E_1$ and $E_2$ are normal, and $M$ is locally the twisted product of the leaves of $E_1$ and $E_2$.
\end{cor}
\begin{proof}
The first results are a direct consequence of the integrability and normality conditions of the previous section.

To see the last piece, recall that $\onab$ describes $H$ completely, and that since $\rP = 0$, $\nabla$ describes $\onab$ completely. Since $\nabla$ is flat, it is entirely described by its preserved sections. Since $\rP = 0$, $\kappa$ must be the torsion of $\nabla$, more specifically the Lie bracket of the flat sections of $T$ fixed by $\nabla$. Now $\kappa(E_j, T) = \kappa(E_j, E_j)$ implies that $\kappa(E_1, E_2) = 0$, demonstrating that the flat sections of $E_1$ fixed by $\nabla$ commute with the flat sections of $E_2$. This makes $\nabla|_{E_j}$ invariant along $E_k$ for $j \neq k$. Moreover, these sections must commute with the flat sections of $B = \{A_1, A_2\}$, meaning that $\nabla$ is invariant along $B$. So each leaf of $E_j$ carries the same free distribution $(M_j, H_j)$. We may divide out by the action of $B$ to get a product manifold $M_1 \times M_2$.

It is then easy to see that $N = H_1 \bigotimes H_2$ caries the same Tractor connection $\onab$ as $M$ does, thus implying they are locally isomorphic.
\end{proof}

\begin{prop}
There exists manifolds with non-flat $\onab$ that fix a minimal co-isotropic subundle $V$, but only if the rank $n$ is $\geq 4$.
\end{prop}
\begin{proof}
If $n = 2,3$, there are no torsion terms, so the flatness of $\nabla$ implies the flatness of $\onab$.

For $n \geq 4$, take the homogeneous model from Section \ref{hom:mod}, and replace $X_1$ and $X_2$ with
\be
X_1' &=& X_1 + x_{12} U_{34} + \frac{1}{2} x_2 U_{21},\\
X_2' &=& X_2 - \frac{1}{2} x_1 U_{12}
\ee
Then $[X_1',X_2'] = U_{12}$, as before, and all the other brackets are as before (since $U_{12}$ only appears in $X_1$ and $X_2$) except for
\be
[U_{12}, X_1'] = U_{34}.
\ee
Then define the flat $\nabla$ by annihilating this new frame for $H$ (hence for $T_{-2}^*$), and $\onab_X = X + \nabla_X$. The only piece of curvature of $\onab$ is
\be
\ora{R}_{U_{12}, X_1'} = Tor^{\nabla}_{U_{12}, X_1'} = U_{34}.
\ee
$\onab$ is normal since the only relevant terms in $\partial^* \ora{R}$ are
\be
\{ U_{12}^*, \ora{R}_{U_{12}, X_1} \} &=& \{ U_{12}^*, U_{23} \} = 0, \\
\{X^*_1, \ora{R}_{X_1, U_{12}} \} &=& \{X_1^*, U_{23} \} = 0, \\
\ora{R}_{ \{U_{12}^*, - \}_{-}, U_{12}} &=& \ora{R}_{0, U_{12}} = 0, \\
\ora{R}_{\{X_1^*, - \}_{-}, X_1}, &&
\ee
and the last term is zero as $\{-,-\}_{-} \subset H$ and $\ora{R}$ is zero on $H \wedge H$. Moreover, if is easy to see that
\be
X_{j} \cdot x_k = \delta_{jk},
\ee
so the hypothesis of Theorem \ref{comb:theo} are fulfilled.
\end{proof}

\section{Fefferman constructions}
Consider a parabolic geometry $(M, \mc{P}, \omega)$ derived from the homogeneous model $G/P$. Assume that there is an inclusion $G \subset \widehat{G}$ with a parabolic inclusion $\wh{P} \subset \wh{G}$ such that $\wh{P} \cap G \subset P$. Assume further that the inclusion $G/(\wh{P} \cap G) \subset \wh{G} / \wh{P}$ is open. Then we may do the \emph{Fefferman construction} on this data. See for example \cite{CR} for details of the original construction.

Define $\wh{M}$ as $\mc{P} / (\wh{P} \cap G)$. The inclusion $(\wh{P} \cap G) \subset \wh{P}$ defines a principal bundle inclusion $i: \mc{P} \hookrightarrow \wh{\mc{P}}$ over $\wh{M}$. Since $\mf{g} \subset \wh{\mf{g}}$, we may extend $\omega$ to a section $\omega'$ of $(T\wh{\mc{P}}^* \otimes \wh{\mf{g}})|_{\mc{P}}$ by requiring that $\omega'(\sigma_A) = A$, for any element $A \in \wh{\mf{g}}$ and $\sigma_A$ the fundamental vector field on $\wh{\mc{P}}$ defined by $A$. We may further extend $\omega'$ to all of $\wh{\mc{P}}$ by $\wh{P}$-equivariance.

Since the inclusion $G/(\wh{P} \cap G) \subset \wh{G} / \wh{P}$ is open, the inclusion $\mf{g} \subset \wh{\mf{g}}$ generates a linear isomorphism $\mf{g}/(\wh{\mf{p}} \cap \mf{g}) \to \wh{\mf{g}} / \wh{\mf{p}}$. At any point $u \in \mc{P}$, $\omega$ is a linear isomorphism $T\mc{P}_u \to \mf{g}$. The previous condition ensures that $\omega'$ is a linear isomorphism $T\wh{\mc{P}}_u \to \wh{\mf{g}}$. This condition extends to all of $\mc{P}$, then to all of $\wh{\mc{P}}$ by equivariance. Consequently $\omega'$ is a Cartan connection.

Dividing out by $P / (\wh{P} \cap G)$ makes $\wh{M}$ into a principal bundle over $M$. It is then easy to see that $\omega'$ is invariant along the vertical vectors of $\wh{M}$ and projects to $\omega$ on $M$. Thus $\omega'$ has holonomy group contained in $G$.

\begin{rem}
The Fefferman construction implies nothing about the relative normalities of $\omega$ and $\omega'$.
\end{rem}
\subsection{Almost-spinorial structures} \label{spin}

There is an evident inclusion of $SO(n+1,n)$ into $SO(n+1,n+1)$. In terms of Dynkin diagrams,
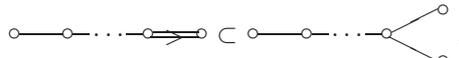
\be
\Balgebra \subset \Dalgebra \ .
\ee

\begin{prop}
There exists a Fefferman construction for this inclusion, where $\wh{M} = M$. In terms of Dynkin diagrams with crossed nodes (see \cite{CartEquiv}), this is
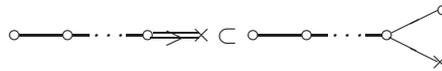
\be
\Bpalgebra \subset \Dpalgebra \,
\ee
and the other parabolic geometry is an almost-spinorial geometry (see \cite{Invar}).

Conversely, any almost spinorial geometry whose Tractor connection preserves a section of the standard Tractor bundle generates a free $n$-distribution on the manifold.
\end{prop}
\begin{proof}
The homogeneous model for the almost-spinorial geometry is $\wh{G} / \wh{P}$ where $\wh{G} = SO(n+1,n+1)$ and $\wh{P}$ is the stabilizer of an isotropic $n+1$ plane. The homogeneous model for a free $n$-distribution are given by $G = SO(n+1,n)$ and $P$ the stabilizer of an isotropic $n$ plane. Since the space $\mbb{R}^{(n+1,n)} \subset \mbb{R}^{(n+1,n+1)}$ must be transverse to every isotropic $n+1$ plane, $\wh{P} \cap G = P$. The open inclusion for the Fefferman condition is equivalent with the statement that $\mf{g}$ and $\wh{\mf{p}}$ are transverse inside $\wh{\mf{g}}$. A simple comparisons of the ranks of all the algebras involved demonstrates that this is the case. This allows us to do the Fefferman construction.

Since $\wh{P} \cap G = P$,
\be
\wh{M} = \wh{\mc{P}} / \wh{P} = \mc{P} / \wh{P} \cap G = \mc{P} / P = M.
\ee
So this almost spinorial structure is on the same manifold as the free $n$-distribution.

Now let $(M, \wh{\mc{P}}, \onab^{as})$ be an almost-spinorial geometry. Let $\mc{T}^{as}$ be its standard tractor bundle. Let $\onab \tau = 0$, for $\tau$ a non-degenerate section of $\mc{T}^{as}$. This gives a reduction of the structure group of $\onab$, making it onto a connection on the principal bundle $\mc{G} \subset \wh{\mc{G}}$ (where these principal bundles have structure groups $G$ and $\wh{G}$ respectively). On $\mc{G}$, $\onab^{as}$ is given by $\omega^{as}$, a one-form with values in $\mf{g}$.

Locally on an open set $U \subset M$, $\mc{G} = G \times U$ and $\wh{\mc{P}} = \wh{P} \times U$. Since $G \cap \wh{P} = P$ for every embedding $G \subset \wh{G}$ given by the preservation of a non-degenerate element,
\be
\mc{G} \cap \wh{\mc{P}} = \mc{P},
\ee
with structure group $P$. Then $\omega^{as}$ restricts to $\mc{P}$, and becomes a Cartan connection on it, inheriting equivariance and point-wise isomorphism.
\end{proof}

The geometric meaning of this is not hard to see. An almost-spinorial geometry has an isomorphism $TM \cong \wedge^2 U$, and there is a projection $\pi: \mc{T}^{as} \to U$. If $\onab^{as} \tau = 0$, this gives us a distribution $H = \pi(\tau) \wedge U \subset TM$. The rank of $U$ must be $n+1$, so the rank of $H$ is $n$; this is our free $n$-distribution.

There is another Fefferman construction that seems relevant here; that given by the inclusion
\be
\Dpalgebrajaw \subset \Bpalgebra \ .
\ee
But except when the first algebra is $D_4$ or $D_3$, parabolic geometries of the the first type are flat if regular and normal (since all their harmonic curvatures have zero homogeneity, see Kostant's solution of the Bott-Borel-Weil theorem \cite{Kostant}). The case of $D_3$ will be dealt with in Section \ref{4:3}.

\section{Free $3$-distributions} \label{4:3}
These are the geometries detailed by Bryant in \cite{bryskew}. They have two properties that distinguish them from the general free $n$-distribution behaviour. First of all, they are torsion free, see Section \ref{har:cur}. Secondly, the almost spinorial structure of Section \ref{spin} is given by $\trione$. However, triality implies that
\be
\trione  \ \ \cong \ \tritwo \ ,
\ee
i.e.~that the almost-spinorial structure is actually a conformal structure, whenever the $SO(4,3)$ structure lifts to a $Spin(4,3)$ structure. This is always true locally.

Paper \cite{bryskew} details the Fefferman construction explicitly. He further shows that if the Tractor connection for the free $3$-distribution is regular and normal, the conformal Tractor connection must be normal as well (regularity is automatic since the conformal parabolic is $|1|$-graded. The holonomy of that conformal Tractor connection must evidently reduce to $Spin(3,4)$.

In fact, the conformal structure is determined by the filtration of $T$ coming from the Tractor connection of the $3$-distribution (see next section). Consequently this local lift globalises for all free $3$-distributions.

\begin{prop}
Conversely, if the normal conformal holonomy of a six manifold $M$ reduces to $Spin(4,3)$, this manifold is the Fefferman construction of a regular normal free $3$-distribution.
\end{prop}
\begin{proof}
Set $\wh{G} = SO(4,4)$, with $\wh{P}$ being $CO(3,3) \rtimes \mbb{R}^{(3,3)}$, the conformal parabolic (defined as the stabiliser of a nul-line in $\mbb{R}^{(4,4)}$. $G = Spin (4,3)$ and $P = GL(3) \rtimes \mbb{R}^{3} \rtimes \wedge^2 \mbb{R}^3$ as before. Let $\onab^c$ and $\omega^c$ be the normal conformal Tractor and Cartan connections.

Let $\wh{\mc{P}}$ be the conformal $\wh{P}$ bundle, $\wh{\mc{P}} \subset \wh{\mc{G}}$ with $\wh{\mc{G}}$ the full structure bundle for $\omega^c$. The holonomy reduction implies that there exists a $G$-bundle $\mc{G} \subset \wh{\mc{G}}$ such that $\omega^c$ reduces to a principal connection on $\mc{G}$.

The action of $Spin(4,3)$ on the nul-lines of $\mbb{R}^{(4,4)}$ is transitive; consequently $G$ and $\wh{P}$ always lie transitively in $\wh{G}$. This means that $\mc{G} \cap \wh{\mc{P}} = \mc{P}$, a $P$-bundle, so $\omega^c$ reduces to a free $3$-distribution Cartan connection -- call it $\omega$.

It remains to show that this Cartan connection is normal. Looking at the homogeneous model, the conformal structure comes from the fact there is a unique conformal class of $\wh{P}$-invariant inner products on $\wh{\mf{g}}/\wh{\mf{p}}$. This implies there is a unique conformal class of $P$ invariant inner products on $\mf{g}/\mf{p}$.

Since $TM = \mc{G} \times_P (\mf{g}/\mf{p})$, this means that the conformal structure on $TM$ depends only on the negative homogeneity components of $\omega$ -- the soldering form, $\omega_-$ (see next section for the geometric details of this).

The curvature $\kappa^c$ of $\omega^c$ can be seen as a $\wh{P}$-invariant map from $\wh{\mc{P}}$ to $\wedge^2 (\wh{\mf{g}}/\wh{\mf{p}})^* \otimes \wh{\mf{g}}$. Similarly, the curvature $\kappa$ of $\omega$ is a $P$-invariant map from $\mc{P}$ to $\wedge^2 (\mf{g}/\mf{p})^* \otimes {\mf{g}}$. On $\mc{P}$, these two curvatures are related by the commuting diagram:
\be
\begin{array}{ccc}
\wedge^2 \mf{g}/\mf{p} &\stackrel{\kappa}{\longrightarrow}& \mf{g} \\
\uparrow & & \downarrow \\
\wedge^2 \wh{\mf{g}}/\wh{\mf{p}} &\stackrel{\kappa^c}{\longrightarrow}& \wh{\mf{g}}.
\end{array}
\ee
Since $\omega^c$ is normal, it is torsion free (see \cite{TCPG} or \cite{mecon}), implying that it maps into $\wh{\mf{p}}$. This means that $\omega$ also maps into $\mf{p}$ -- so is also torsion-free. This means that $\kappa$ is of homogeneity $\geq 2$, consequently -- since $\partial^*$ respects homogeneity -- $\partial^* \kappa$ is of homogeneity $\geq 2$.

Now, by \cite{capslo}, any Cartan connection $\omega$ with curvature $\kappa$ such that $\partial^* \kappa$ is of homogeneity $\geq l \geq 0$ differs from the normal Cartan connection $\omega'$ by a section $\Phi \in \Omega^1(\mc{P},\mf{g})$ of homogeneity $\geq l$.

So here we have $\omega + \Phi = \omega'$, with $\omega'$ normal and $\Phi$ of homogeneity $\geq 2$. This means that $\omega'$ and $\omega$ have the same soldering form (as the soldering form is of strictly negative homogeneity), thus that the conformal structure that they both generate are the same. Since the conformal Fefferman construction for $(\mc{P}, \omega')$ must be normal (since $\omega'$ is), it must \emph{be} $(\wh{\mc{P}}, \omega^c)$. This means that $\Phi = 0$, hence that $\omega$ is normal.
\end{proof}

\subsection{Geometric equivalence}

Given a free $3$-distribution on manifold $M$, the conformal structure can be recovered directly from the decomposition of $T \cong T_{-2} \oplus H$ given by any Weyl structure. Let $\sigma$ be any local never-zero section of $\wedge^3 H$. Then there is a map $g: T_{-2} \otimes H \to \wedge^3 H$ given by the isomorphism $T_{-2} \cong \wedge^2 H$. Extend $g$ to a section of $(\odot^2 T^*) \otimes \wedge^3 H$ by the inclusion $T_{2}^* \otimes H^* \subset \odot^2 T^*$. Then $g\sigma^{-1}$ is a metric on $M$. This depends on the choice of the section $\sigma$, so actually defines a conformal structure. It is then easy to see that $g$ is invariant under the action of a one-form $\Upsilon$, (as $g(U+Y,X) = g(U,X) + g(X,Y) = g(U,X)$ for any sections $X$ and $Y$ of $H$ and any section $U$ of $T_{-2}$). So this conformal structure does not depend on the choice of Weyl structure, only on the filtration of $T$ (which depends on the Cartan connection).

The algebra $\mf{spin}(4,3) \subset \mf{so}(4,4)$ is defined by preserving a generic four-form $\lambda$ on $V = \mbb{R}^{(4,4)}$, see \cite{baumspin}. Let $\mc{T}^C$ be the standard conformal Tractor bundle on $M$ (see \cite{TBIPG} or \cite{mecon} for more details on conformal geometries). If the conformal Tractor connection $\onab^C$ has holonomy algebra reducing to $\mf{spin}(4,3)$, then there exists a generic four-form $\nu \in \Gamma(\wedge^4 \mc{T}^C)$ such that
\be
\onab^C \nu = 0.
\ee
There is a natural projection on $\mc{T}_C$, coming from its filtration
\be
\mc{T}_C  \to \mc{E}[1] \oplus T[-1] \to \mc{E}[1].
\ee
Here $\mc{E}[1]$ is a density bundle, $\mc{E}[\alpha] = (\wedge^6 T^*)^{\frac{\alpha}{-6}}$, and $T[-1] =T \otimes \mc{E}[-1]$. This implies that there is a well defined projection $\pi: \wedge^4{\mc{T}_C} \to (\wedge^3 T)[-2]$. It turns out that $\pi(\nu)$ is decomposable, and so defines a distribution $H^*$ of rank three in $T^*[2/3]$. Since a distribution is unchanged by a change of scale, this is actually a distribution in $T^*$, with dual distribution $H \subset T$. This $H$ is precisely that defining the Bryant structure; the maximal non-integrability derives from the properties of $\nu$ and $\onab^C$.

\subsection{$G_2'$ structures}
The most natural subgroup of $Spin(4,3) \subset SO(4,4)$ is $G_2'$, defined as the subgroup of $Spin(3,4)$ that preserves a non-isotropic element $e$ in $\mbb{R}^{(4,4)}$.

There are many equivalent definitions. For examples, $G_2'$ is equivalently described as the subgroup of $SO(4,3)$ that preserves a generic three-form $\theta$ on $\mbb{R}^{(4,3)}$; $\theta$ is just $e \llcorner \lambda$. Alternatively, it is the automorphism group of $\mbb{O}'$, the split Octonions. It acts irreducibly on the seven dimensional space $V = Im \mbb{O}'$. The split Octonions carry a natural inner product $N$, generated from the norm $N(x,x) = x \overline{x}$. This quadratic form is multiplicative, and is of signature $(4,4)$. The identity element $1 \in \mbb{O}'$ is of positive norm squared, and is orthogonal to $V$; thus $V$ is of signature $(3,4)$. This algebra is alternative; this means that the alternator
\be
[x,y,z] = (xy)z - x(yz)
\ee
is totally anti-symmetric in its three entries. We can use $N$ to make $[,,,]$ into an element of $\wedge^3 V^* \otimes V^*$; it turns out to be skew in all four entries, and equal to $* \theta$ where $*$ is the Hodge star generated by $N$.

The three-form $\theta$ itself is given by
\be
\theta(x,y,z) = N(xy,z).
\ee
The properties of the split Octonions force this to be skew in all three arguments.

Now assume that our free $3$-distribution has a normal Tractor connection $\onab$ with a holonomy reduction to $G_2'$. By the conformal Fefferman construction, the conformal structure will be given by a manifold that is conformally Einstein and whose metric cone carries a $G_2'$ structure (see \cite{meein}). Such manifolds do exist -- for instance, $SL(3, \mbb{R}) /T^2$ where $T^2$ is a maximal torus, is one example \cite{brymetric}. Here, the free $3$-distribution would be chosen at $Id \in SL(3, \mbb{R})$ as the span of
\be
H_{Id} = \left( \begin{array}{ccc} 0 & a &0 \\ 0&0&b \\ c&0&0 \end{array} \right),
\ee
and extended to the whole manifold by Lie multiplication. Note that $\{ H_{Id}, \mf{T}^2 \} \subset H_{Id}$, for $\mf{T}^2$ the tangent space to the maximal torus, so this extension is well defined.

It is not unique, however. We could have used the transpose of $H_{Id}$ instead. Note that $H_{Id}^t = \{H_{Id}, H_{Id} \}$. This will be an important property of $G_2'$ structures on a free $3$-distribution.

\begin{prop}
There are three orbits of isotropic $3$-planes in $\mbb{R}^{(4,3)}$ under the action of $G_2'$. -- two open, one closed. The closed orbit is distinguished by the fact that $\theta(x,y,z) = 0$ for all elements in an isotropic $3$-plane inside this orbit.
\end{prop}
\begin{proof}
Let $B \subset Im \mbb{O}'$ be an isotropic $3$-plane. For any element $x$ of $B$, $0 = N(x,x) = x \times \overline{x} = - x \times x$. This implies that for any elements $x$ and $y$ of $B$,
\be
0 &=& 2 N(x,y) = N(x+y,x+y) - N(x,x) - N(y,y) \\
&=& (x+y) \times (\overline{x} + \overline{y}) - x \times \overline{x} - y \times \overline{y} \\
&=& -(x+y) \times (x+y) \\
&=& -x \times y - y \times x.
\ee
Thus an isotropic $3$-plane is defines as a subset of $Im \mbb{O}'$ where every element squares to zero, and anti-commute. There are two situations to be covered:

\begin{enumerate}
\item There exists a basis $\{x,y,z\}$ for $B$ such that $\lambda(x,y,z) = 1$.
\end{enumerate}
The set of all such $B$ is evidently open in the set of all isotropic $3$-planes. We aim to show $G_2'$ is transitive on this set.

\begin{lemm}
The elements span of $x$, $y$ and $z$ under split Octonionic multiplication generate all of $\mbb{O}'$.
\end{lemm}
\begin{lproof}
Since the split Octonions are alternative, the multiplicative span of any two elements is associative. Hence
\be
(xy)(xy) = x(xy)y = -x(xy)y = -(xx)(yy) = 0.
\ee
This is true for any elements $x$, $y$ in B. Thus $C = B \times B$ is isotropic, so of maximum dimension three. The relation
\be
1 = \lambda(x,y,z) = N(xy,z),
\ee
implies that $xy$ is orthogonal to $x$, $y$, but not to $z$. We may cyclically permute $x$, $y$ and $z$ here, thus demonstrating that $C$ is of dimension three.

In fact $N(x \wedge y \wedge z, yz \wedge zx \wedge xy) = 1$, so
\be
\lambda(yz,zx,xy) = -1.
\ee

Then define $a = (xy)z - z(xy)$. Now $\overline{a} = \overline{z}(\overline{xy}) - (\overline{xy}) \overline{z} = z(xy) - (xy)z = -a$, so $a$ is purely imaginary. We make the claim that $x,y,z,xy,zx,yz$ and $a$ span $Im \mbb{O}'$ and that the split Octonion multiplication of these elements is completely determined.

First of all, the squares of $x,y,z,xy,zx$ and $yz$ are all zero, as are all the terms $x \times xy$, $x \times zx$, $y \times yz$, $y \times xy$, $z \times zx$ and $z \times yz$. Now consider $b = (xy) \times (yz)$. How $N(b,x) = \lambda(xy,yz,y)= N((yz)y,xy) = N(0,xy) = 0$. Similarly $N(b,z) = N(b,y) = 0$; thus $b \in B^{\perp}$. Moreover $N(b,b) = N(xy,xy)N(yz,yz) = 0$, so $b \in B$. Trying to extract the $x$, $y$ and $z$ components of $b$, we find $N(b,xy) = 0$, $N(b,zx) = -1 $ $N(b,yz) = 0$, so $b = -y$. Similar reasoning demonstrates that
\be
yz \times zx = -z \ \ \textrm{and} \ \ zx \times xy = -x.
\ee
Similar manipulations, using the associator $*\lambda$, show that
\be
xa = x = -ax, & ya = y =-ay, & za = z = -za \\
a(yz) = xy = -(yz)a, & a(zx) = zx =-(zx)a, & a(xy) = xy = -(xy)a
\ee
and
\be
a = (yz)x - x(yz) = (zx)y - y(zx).
\ee
And finally $a \times a = 1$.

It often helps to work with an explicit description of split Octonion multiplication. Here is one, due to Zorn. Here, a split Octonion is represented by the ``matrix''
\be
x= \left( \begin{array}{cc} a & \bf{v} \\ \bf{w} & b \end{array}\right)
\ee
with $a$ and $b$ real numbers and $\bf{v}$ and $\bf{w}$ vectors in $\mbb{R}^3$. The norm squared $N(x,x)$ is the ``determinant'' $ab - \bf{v} \cdot \bf{w}$. Multiplication is given by
\be
\left( \begin{array}{cc} a & \bf{v} \\ \bf{w} & b \end{array}\right) \times \left( \begin{array}{cc} a' & \bf{v}' \\ \bf{w}' & b' \end{array}\right) = \left( \begin{array}{cc} aa' + \bf{v}\cdot\bf{w}' & a\bf{v}' + b'\bf{v} + \bf{w} \wedge \bf{w}' \\ a'\bf{w} + a \bf{w}' - \bf{v} \wedge \bf{v}' & bb' +\bf{v}' \cdot \bf{w} \end{array}\right).
\ee
With $\cdot$ and $\wedge$ the ordinary dot and cross products on $\mbb{R}^3$. The imaginary split Octonions are those where $a  = -b$.
\end{lproof}

So if we call $\{x,y,z\}$ an \emph{Octonionic triple}, then an element $g$ of $G_2'$ is entirely determined by $\{g(x),g(y),g(z)\}$. Conversely, for any two Octonionic triple, there is an element of $G_2'$ mapping one to the other.

Moreover, if $G_{B} \subset SL(7,\mbb{R})$ is the stabiliser of $B$, $G_2' \cap G_B$ is the permutation group of the Octonionic triples in $B$ -- consequently $G_2' \cap G_B = SL(3,\mbb{R})$.

\begin{enumerate}
\setcounter{enumi}{1}
\item For all $x,y,z \in B$, $\lambda(x,y,z) = 0$.
\end{enumerate}
The set of all such $B$ is closed in the set of all isotropic $3$-planes, complementary to the previous orbit, and with empty interior. We aim to show $G_2'$ is transitive on this set.

As in the previous examples, $(xy)(xy) = 0$. And $C = B \times B$ is isotropic. However
\be
0 = \lambda(x,y,z) = N(xy,z),
\ee
implying that $C \subset B^{\perp}$. Since $C$ is isotropic, $C \subset B$. An inspection of the explicit form of split Octonion multiplication demonstrates that there does not exist a three plane on which $\times$ is totally degenerate. So $C \neq 0$. Let $z \in C$. Now $z = xy$ for elements $x$ and $y$ in $B$. Since elements of $B$ square to zero, $x \neq y$. Since the multiplicative span of any two elements is associative, $z \neq y$ and $z\neq $. Furthermore, $z$ cannot be in the span of $x$ and $y$, since $x (r_1 x + r_2 z) = r_1 xx + r_2 x(xy) = 0$ for all real $r_j$. So $x$, $y$ and $z$ form a basis for $B$, and the relations
\be
xy = z, xz =0, yz = 0, xx = yy = zz = 0,
\ee
determine multiplication on $B$ completely. In fact, $B$ is determined by $z$. This can be seen from the fact that $G_2'$ is transitive on the set of isotropic element of $Im \mbb{O}'$, so we may set
\be
z = \left( \begin{array}{cc} 0 & e_1 \\ 0 & 0 \end{array} \right),
\ee
where $e_1, e_2, e_3$ is a basis for $\mbb{R}^3$. Then the two sided kernel of the multiplications $\times z, z \times : Im \mbb{O}' \to \mbb{O}'$ is spanned by
\be
z, \left( \begin{array}{cc} 0 & 0 \\ e_2 & 0 \end{array} \right), \left( \begin{array}{cc} 0 & 0 \\ e_3 & 0 \end{array} \right).
\ee
Since $B$ is in the two-sided kernel of multiplication by $z$, and is isotropic, it must be precisely the span of these elements. Since $B$ is determined by $z$, and since $G_2'$ is transitive on isotropic elements of $Im \mbb{O}'$, $G_2'$ must be transitive on the set of isotropic $3$-planes $B$ on which $\lambda$ vanishes.

\end{proof}

\begin{theo}
Let $M$ be a free $3$-distribution manifold with normal Tractor connection $\onab$, with the holonomy group of $\onab$ reducing to $G_2'$. Then, on an open, dense set of $M$, there is a unique Weyl structure $\nabla$ defined by this information. This Weyl structure determines a splitting of $T = T_{-2} \oplus H$. Then $H'=T_{-2}$ is a free $3$-distribution. And the normal Tractor connection determined by $H'$ is isomorphic to $\onab$.
\end{theo}
\begin{proof}
If $\onab$ has holonomy contained in $G_2'$, then it comes from a connection an a principle $G_2'$-bundle $\mc{G}_2'$. Let $\mc{A}' = \mc{G}_2' \times_{G_2'} \mf{g}_2'$. The inclusion $G_2' \subset SO(3,4)$ generates inclusions $\mc{G}_2' \subset \mc{G}$ and $\mf{g}_2' \subset \mf{g}$, thus an inclusion $\mc{A}' \subset \mc{A}$. And by definition $\onab$ preserves $\mc{A}'$, and $L$, a three-form on $\mc{T}$.

Since $\onab$ has holonomy contained in $G_2$, it also preserves split Octonionic multiplication on $\mc{T}$. Designate this multiplication by $\times$. By the previous proposition, on an open, dense subset of $M$, the canonical $H^* \subset \mc{T}$ generates all of $\mc{T}$ by $\times$. We have a well defined subundle of $\mc{T}$, $K = H^* \times H^*$. Since $K$ and $H^*$ are transverse, the projection $\pi^2$ maps $K$ isomorphically to $H$. Then let $\nabla$ be the (unique) strongly $K$-preferred connection. In the splitting it defines, set $H' = T_{-2}$.

Now consider $\mc{A}'_0 \subset \mc{A}'$, the subundle of $\mc{A}'$ that stabilises $H^*$ (and $K$). This must be a $\mf{sl}(3,\mbb{R})$ bundle, since the subgroup of $G_2'$ that preserves a generic isotropic $3$-plane is $SL(3,\mbb{R})$. By the way we have chosen our current splitting, $\mc{A}'_0 \subset{A}_0$. Consequently $\nabla$ preserves a volume form, and thus $H' \cong H \wedge H \cong H^*$. Thus under the action of $\mc{A}'_0$,
\be
\mc{A} = H' \oplus H \oplus \mc{A}_0' \oplus \mbb{R} \oplus H' \oplus H.
\ee
Since $\mc{A}'$ is of rank $14$, since $G_2'$ is fourteen dimensional, there are three possibilities for the structure of $\mc{A}'$
\be
\mc{A}' &=& \mc{A}_0' \oplus H \oplus H' \\
\mc{A}' &=& \mc{A}_0' \oplus H \oplus H \\
\mc{A}' &=& \mc{A}_0' \oplus H' \oplus H'.
\ee
But the last two possibilities are not algebraically closed, so $\mc{A}'$ must be of the first type. It is also simple, which means that it cannot be of pure positive or negative homogeneity. A bit of experimentation then shows that the only possibility for $\mc{A}'$ is that it is composed of elements of the form
\be
(X',X,\Theta, X', X),
\ee
for $X \in \Gamma(H)$, $X' \in \Gamma(H')$, $\Theta \in \Gamma(\mc{A}'_0)$ (the same could instead be deduced from the split Octonionic multiplication in this manifold). Since $\onab$ must preserve this bundle, and that $\nabla$ already does, $\rP_{12}$ must be the identity on $H$, $\rP_{21}$ the identity on $H'$ and $\rP_{11} = \rP_{22} = 0$. Notice that $\nabla \rP = 0$ -- this is similar to, though not identical to, an Einstein involution \cite{meein}.

This implies that the Cartan connection $\omega$ decomposes as
\be
\omega = \omega_{-2} + \omega_{-1} + \omega_0 + \omega_1 + \omega_2,
\ee
with $\omega_{-2} = \omega_1$ and $\omega_{-1} = \omega_2$.

Since $\onab$ is torsion free, $\kappa(T \wedge T)$ must take values in $\mc{A}_{(0)} \cap \mc{A}' = \mc{A}'_{0}$. The harmonic curvature component of $\kappa$ (see Section \ref{har:cur}) is in $H \otimes H' \otimes \mc{A}'_0$. The only other possible curvature component of $\kappa$ is the higher homogeneity $H' \wedge H' \otimes \mc{A}'_0$. Since $\rP_{22} = 0$, this is precisely the $R^{\nabla}_{22}$, where $R^{\nabla}$ is the curvature of $\nabla$.

\begin{lemm}
$\kappa_{22} = R^{\nabla}_{22} = 0$.
\end{lemm}
\begin{lproof}
Designate $\kappa_{22}$ by $\kappa'$. The Bianchi identity for $\onab$ is $d^{\onab} \kappa = 0$, where $d^{\onab}$ is $\onab$ on $\mc{A}$ twisted with the exterior derivative $d$ on $\wedge^2 T^*$. For $X'$ and $Y'$ sections of $H'$ and with $Z$ a section of $H$,
\be
0 = (d^{\onab} \kappa)_{X',Y',Z}&=& (d^{\nabla} R^{\nabla})_{X',Y',Z} + \{\kappa,-\}_{X',Y',Z} + \{\kappa,\rP(-)\}_{X',Y',Z} \\ 
&=& \{\kappa_{12},-\}_{X',Y',Z} + \{\kappa_{12},\rP(-)\}_{X',Y',Z} + \kappa_{X',Y'}'  \cdot Z + \kappa_{X',Y'}'  \cdot \rP(Z)
\ee
Now the expression $\kappa_{X',Y'}'  \cdot \rP(Z)$ is the only component taking values in $\mc{A}_{2}$, so it must vanish. This implies that $\kappa' = 0$.
\end{lproof}

Now we have $\kappa$ as a section of $H \otimes H' \otimes \mc{A}'_{(0)}$. In particular $\kappa(H \wedge H) = 0$. Recall the definition of normality; that $\partial^* \kappa = 0$, where
\be
(\partial^* \kappa)(X) = \sum_{l} \{ Z^l, \kappa_{(Z_l, X)} \} - \frac{1}{2} \kappa_{(\{Z^l,X \}_{-}, Z_l)},
\ee
for $(Z_l)$ a frame for $T$ and $(Z^l)$ a dual frame for $T^*$. Now $\{Z^l,X \}_{-} \wedge Z_l$ is zero or a section of $H \wedge H$ for all $X$ and $Z_l$. Thus the normality of $\kappa$ is entirely encapsulated in
\be
\sum_{l} \{ Z^l, \kappa_{(Z_l, X)} \},
\ee
or, in other words, in the fact that $\kappa$ is trace free.

Now consider $\onab$ as a principal connection on $\mc{G}$, forgetting about the inclusion $\mc{P} \subset \mc{G}$. We may define an alternative inclusion $\mc{P}' \subset \mc{G}$ by using $K$ as the canonical subundle of $\mc{T}$. By our previous results, the new soldering form is now $\omega_{-1} + \omega_{-2}$ rather than $\omega_{-2} + \omega_{-1}$ -- and this is a proper soldering form, meaning that $\onab$ is a Tractor connection for the distribution $H'$. The curvature of $\onab$ is still $\kappa$, though the new soldering form sends $H \otimes H'$ to $H' \otimes H$. Under this new identification, $\kappa$ thus remains a trace-free section of $H' \otimes H \otimes \mf{A}'_{(0)}$. Thus if $\partial^{*'}$ is the operator for the new parabolic,
\be
\partial^{*'} \kappa = \textrm{Trace } \kappa = 0.
\ee
Thus $\onab$ is normal as the Tractor connection generated by $H'$.
\end{proof}

\subsection{CR structures}
We aim to show here that there is a Fefferman construction for $\wh{G} = SO(4,3)$, $\wh{P}$ stabilises an isotropic $3$-plane, and $G = SO(4,2)$ while $P = (SO(2) \oplus GL(2)) \rtimes (\mbb{R}^2 \otimes \mbb{R}^{(2)}) \rtimes \wedge^2 \mbb{R}^2$ stabilises an isotropic $2$-plane.

Let $V = \mbb{R}^{(4,3)}$ and $B \subset W$ be an isotropic $3$-plane whose inclusion defines $\wh{P} \subset \wh{G}$. Let $W \cong \mbb{R}^{(4,2)}$, and fix an inclusion $W \subset V$ that defines $G \subset \wh{G}$.

Because of their signatures, $W$ and $B$ must be transverse, so their intersection $C = W \cap B$ is an isotropic $2$-plane. Defining $P$ as the stabiliser of $C$, it is evident that $G \cap \wh{P} \subset P$.

Now let $B'$ be the orthogonal projection of $B$ onto $W$. By construction, $C \subset B' \subset C^{\perp}$. The bundle $B'$ is equivalently defined by a line through the origin in $C^{\perp} / C$. The group $P$ acts via $SO(2)$ on this space of lines. Thus $G \cap \wh{P}$ lies as a codimension one subgroup in $P$. Then $\wh{G}$ is of dimension $21$, $\wh{P}$ of dimension $15$, $G$ also of dimension $15$, $P$ of dimension $10$ and $G \cap \wh{P}$ of dimension $9$. This implies that $G$ and $\wh{P}$ are transverse in $\wh{G}$, hence that the inclusion $\mf{g}/(\mf{g} \cap \wh{\mf{p}} \to \wh{\mf{g}} / \wh{\mf{p}}$ is open. Thus we may do the Fefferman construction.

\begin{defi}[CR] A CR manifold is given by a contact distribution $K \subset TN$ with a complex structure $J$ on $K$. If $Q = TN/K$, and $q: TN \to Q$ is the obvious projection, there is a skew symmetric map $\mc{L}: K \times K \to Q$ given by $\mc{L}(X,Y) = q([X,Y])$ where $X$ and $Y$ are sections of $K$.

Integrability comes from using $J$ to split $K \otimes \mbb{C}$ as $K^{1,0} \oplus K^{0,1}$; the CR structure is integrable if $K^{0,1}$ is closed under the Lie bracket. This implies that $\mc{L}$ is of type $(1,1)$, that is that $\mc{L}(JX,JY) = \mc{L}(X,Y)$.
\end{defi}

\begin{theo}
The geometries modelled on $G/P$ are the 5 dimensional split signature CR geometries. If the CR structure is integrable and the Cartan connection is normal, the Cartan connection on the free $3$-distribution coming from the Fefferman construction is also normal.

Conversely, if the holonomy group of a normal Cartan connection for a free $3$-distribution reduces to $SO(4,2)$, it is the Fefferman construction over an integrable split signature CR geometry with normal Cartan connection.
\end{theo}
The rest of this section will be devoted to proving that theorem.

The first statement -- that the $G/P$ geometries are the CR geometries -- from the fact that the representation of $P$ as a parabolic is $\CR$, the same as for CR structures, combined with the following lemma:
\begin{lemm}
$Spin(4,2)_0 = SU(2,2)$.
\end{lemm}
\begin{lproof}
Consider the action of $SU(2,2)$ on $V = \mbb{C}^{(2,2)} \wedge \overline{\mbb{C}}^{(2,2)}$. $V$ carries a real structure on it from the action of the K{\"a}hler form, and a natural $(4,2)$ signature metric. Since $SU(2,2)$ is simple, and its action is non-trivial on this space, there is an inclusion
\be
\mf{su}(2,2) \hookrightarrow \mf{so}(4,2).
\ee
And then dimensional considerations imply that this is an equality.

The maximal compact subgroup of $SO(4,2)_0$ is $S(O(4) \times O(2))_0$; the maximal compact subgroup of $SU(2,2)$ is $S(U(2) \times U(2))$. This means that the fundamental groups of the Lie groups are:
\be
\pi_1(SO(4,2)_0) &=& \mbb{Z}_2 \oplus \mbb{Z} \\
\pi_1(SU(2,2)) &=& \mbb{Z}. \\
\ee
Consequently $Spin(4,2)_0 = SU(2,2)$.
\end{lproof}
Then it is easy to see that $P$ is the stabiliser of a complex nul-line in $\mbb{C}^{(2,2)}$, demonstrating that these are split signature CR geometries (see \cite{CR}).

In order to demonstrate the normality conditions, we shall use both this Fefferman construction and the conformal Fefferman construction. Let $\wh{\wh{G}} = SO(4,4)$, with $\wh{\wh{P}}$ the stabiliser of a nul line.  In details, if $(\mc{P}, \omega)$ is a split signature CR geometry, we have three related structures:
\be
(\ \mc{P}\ , \ \ \omega\ ) \ \ \ \ \ (\ \wh{\mc{P}}\ , \ \ \wh{\omega}\ ) \  \ \ \  \ (\ \wh{\wh{\mc{P}}}\ , \ \ \wh{\wh{\omega}}\ ),
\ee
where $\wh{\omega}$ is the Tractor connection for a free $3$-distribution while $\wh{\wh{\omega}}$ is a conformal Tractor connection.

We know that ${\wh{\omega}}$ is normal if and only if $\wh{\wh{\omega}}$ is normal. But now consider the total Tractor connection, from $G$ to $\wh{\wh{G}}$. This is determined by how $G$ lies in the larger group.

If we complexify everything, we have $Spin(6) \subset Spin(7) \subset SO(8, \mbb{C})$. The spin representations of $Spin(6)$ are isomorphic with $\mbb{C}^4$, so decompose $\mbb{C}^8$ into two distinct components.

This implies that the action of $Spin(4,2) \subset Spin(4,3) \subset SO(4,4)$ on $\mbb{R}^{(4,4)}$ either decomposes it into two four dimensional components, or is irreducible on it (and preserves a complex structure on it). However, $SU(2,2) = Spin (4,2)$ does not have any four dimensional real representations (apart from the trivial one). Consequently the inclusion $SU(2,2) \subset SO(4,4)$ is the standard inclusion.

This means that the inclusion $G \subset \wh{\wh{G}}$ is the standard one. This generates $\wh{\wh{\omega}}$ via the Fefferman construction. But this Fefferman construction has to be the standard one. This implies that $\wh{\wh{\omega}}$ is normal if and only if $\omega$ is normal and the CR structure is integrable (see \cite{CR}, \cite{leitnersu} and \cite{leitnercom}). Consequently, $\wh{\omega}$ is normal if and only if $\omega$ is normal and the CR structure is integrable.

\begin{rem}
The inclusion $SU(2,2) \subset Spin(4,3)$ can be seen directly. $Spin (4,3)$ is defined as preserving a generic four-form $\lambda$ on $\mbb{R}^{(4,4)}$ (see \cite{baumspin}). $SU(2,2)$ on the other hand, preserves a K{\"a}hler form $\mu$, which can be seen as a section of $\wedge^2 V$ that is conjugate linear with respect to the volume form. It also preserve a complex volume form $v \in \Gamma(\wedge^{(4,0)} V_{\mbb{C}})$. The inclusion of $SU(2,2)$ into $Spin(3,4)$ is given by the generic four form:
\beqa \label{form:su}
Re(v) - (\mu)^2.
\eeqa
\end{rem}

Now we need to show the converse. Let $(M, \wh{\mc{P}}, \onab)$ be a normal Cartan connection for a free $3$-distribution. Assume the holonomy group of $\onab$ reduces to $SO(4,2)$ -- equivalently, that there is a section $\tau$ of $\mc{T}$, of negative norm squared, such that $\onab \tau = 0$. Define $R \in \Gamma(H)$ as $\pi^2(\tau)$; since $\tau$ is of negative norm squared, $R$ is never-zero.

Define the bundle $\wt{K}$ as $H \oplus [H,R]$. It is a bundle of rank five. Let $N$ be the manifold got from $M$ by projecting along the flow of $R$.
\begin{prop}
$N$ carries a CR structure, and the contact distribution $K$ in $TN$ is the projection of $\wt{K}$ to $N$. The complex structure $J$ on $K$ is given by the action of $R$.
\end{prop}

\begin{proof}
We first need to show that $[\wt{K}, R] = \wt{K}$. Let $L$ be the line subundle of $\mc{T}$ generated by $\tau$. Pick any $L$-preferred connection $\nabla$. This obeys the properties of Theorem \ref{V:pref} -- $\tau = (\alpha,0,R)$, $\alpha(R) = 1$ and $\nabla_X \alpha = \nabla_X R = 0$ for any section $X$ of $H$ while $\nabla_U R = - \{U,\alpha\}$ for $U$ a section of $T_{-2}$.

We may choose sections $X$ and $Y$ of $H$ that obey the following properties:
\begin{enumerate}
\item $X$, $Y$ and $R$ form a frame of $H$,
\item $\nabla_R X = \nabla_R Y = 0$,
\item $\alpha(X) = \alpha(Y) = 0$, 
\end{enumerate}
(for instance, we could define $X$ and $Y$ obeying the algebraic properties on a submanifold transverse to $R$, and extend by parallel transport along $R$; then the relation $\nabla_R R = \nabla_R \alpha = 0$ ensures the algebraic properties are preserved). Since $\onab$ is torsion-free,
\be
[R, X] &=& \nabla_R X - \nabla_X R - \{R,X\} \\
&=& \{R,X\}.
\ee
Similarly,
\be
[R,[R,X]] &=& \nabla_R \{R,X\} - \nabla_{\{R,X\}} R - 0 \\
&=&  \{\{R,X\}, \alpha \} = -X.
\ee
The same hold for $Y$ and $\{R,Y\}$. Consequently $[R, \wt{K}] = \wt{K}$ and $\wt{K}$ projects to a distribution $K$ in $N = M/R$. This distribution must be a contact distribution, by the properties of the Lie bracket on $\wt{K}$. Let $r$ be any coordinate on $M$ such that $R\cdot r = 1$. Then the vector fields
\beqa \label{sec:CR}
\cos(r)X - \sin(r)\{R,X\},& \sin(r)X - \cos(r)\{R,X\} \\
\cos(r)Y - \sin(r)\{R,Y\},& \sin(r)Y - \cos(r)\{R,Y\}
\eeqa
are $R$-invariant, hence lifts of vector fields in $K$. Since we have this explicit form, we can see that the Lie bracket of $X$ and $Y$ with $R$ generates an endomorphism of $\wt{K}$ that descends to an automorphism $J$ of $K$, squaring to minus the identity.

Changing to another $L$-preferred connection will change $X$ and $Y$ by adding multiples of $R$. This will change neither their projections nor the properties of $J$. This means that the CR structure is well defined.

And $M$ must be the Fefferman construction over this CR structure. This implies that CR structure must be integrable and that $\onab$ descends to a normal CR Tractor connection on $N$.
\end{proof}

One interesting consideration: though the free $3$-distribution defines the CR structure uniquely, the converse is only true up to isomorphism. Any diffeomorphism $\phi: M \to M$ generated by a flow on $R$ will change the distribution $H$, but since $\phi$ projects to the identity on $N$, it leaves the underlying CR structure invariant.

The forgoing means that all the results on CR holonomy (equivalently, conformal holonomy contained in $\mf{su}(2,2)$) have equivalent formulations in terms of free $3$-distributions. See papers \cite{leitnersu}, \cite{CR} and \cite{leitnercom}; paper \cite{mecon} has some Einstein examples as well. This implies, for instance, that holonomy reductions to $SU(2,1)$ exist (whenever $N$ is a Sasaki-Einstein manifold with the correct signature and sign of the Einstein coefficient). From the free $3$-distribution point of view, this corresponds to a complex structure on the complement of $\tau$ in $\mc{T}$.

Similar consideration exist for a holonomy reduction to $Sp(1,1) \subset SU(2,2)$, with the quaternionic analogue of CR spaces.

\subsection{Lagrangian contact structures}
Lagrangian contact structures (see for example \cite{complexconf}) geometries generated by another real form of the parabolic that models CR structures.

\begin{defi}
A Lagrangian contact structure is given by a contact distribution $K$ on a manifold of dimension $2m+1$, together with two bundles $E$ and $F$ of rank $m$ such that $K = E \oplus F$, $[E,E] \subset K$ and $[F,F] \subset K$.

The structure is integrable if both $E$ and $F$ are integrable. The parabolics are given by $G = SL(m, \mbb{R})$ while $P = (\mbb{R} \oplus GL(m-2, \mbb{R}) \rtimes (\mbb{R}^{m} \oplus \mbb{R}^{m*}) \rtimes \mbb{R}$.
\end{defi}

Then we have:
\begin{lemm}
$Spin(3,3)_0 = SL(4,\mbb{R})$
\end{lemm}
\begin{lproof}
Consider the action of $SL(4,\mbb{R})$ on $V = \wedge^2 R^4$. Since $SL(4,\mbb{R})$ preserves a volume form which is an element of $\wedge^4 \mbb{R}^{4*} \subset \odot^2 (\wedge^2 \mbb{R}^4)^*$, it preserves a metric on $V$, of split signature. Then since $SL(4,\mbb{R})$ is simple and acts non-trivially, we get an algebra inclusion $\mf{sl}(4,\mbb{R}) \subset \mf{so}(3,3)$ and the dimensions imply equality.

The maximum compact subgroup of $SL(4,\mbb{R})$ is $SO(4)$ while the maximum compact subgroup of $SO(3,3)_0$ is $S(O(3) \times O(3))_0$. Consequently
\be
\pi_1 (SL(4,\mbb{R})) = \mbb{Z}_2, \\
\pi_1 (SO(3,3)_0) = \mbb{Z}_2 \times \mbb{Z}_2,
\ee
demonstrating the result.
\end{lproof}

Given this, the results for CR structures go through almost verbatim to this new setting. There is one subtlety, however: $\mbb{R}^{(3,3)}$ need not be transverse to a given isotropic $3$-plane in $\mbb{R}^{(4,3)}$. So we may often need to restrict our results to open dense subsets of our manifolds. Also the inclusion $GL(4,\mbb{R}) \subset Spin(4,3) \subset SO(4,4)$ is now the standard inclusion that decomposes $\mbb{R}^{(4,4)}$ as $\mbb{R}^4 \oplus \mbb{R}^{4*}$.

Summarising all the results:
\begin{theo}
Let $N$ be a five dimensional integrable Lagrangian contact manifold. Then there is a Fefferman construction for $N$ to a free $3$-distribution on a manifold $M$. The Tractor connection on $M$ is normal if and only if the Tractor connection on $N$ is normal.

Conversely, if $M$ is a free $3$-distribution geometry with normal Tractor connection $\onab$, and the holonomy group of $\onab$ reduces to $SO(3,3)$, then \emph{an open dense set of $M$} is the Fefferman space of a integrable, normal Lagrangian contact manifold.
\end{theo}

\bibliographystyle{alpha}
\bibliography{ref}

\end{document}